\documentclass[a4paper,fleqn]{cas-sc}

\usepackage[authoryear,longnamesfirst]{natbib}
\usepackage[utf8]{inputenc}
\usepackage{amsmath}
\usepackage{amsthm}
\usepackage{longtable}
\usepackage{booktabs}
\usepackage[algoruled,linesnumbered]{algorithm2e}
\usepackage{microtype}
\usepackage{tikz}

\newcommand{\BPP}{{BPP}}
\newcommand{\BPPS}{{BPPS}}
\newcommand{\TP}[1]{\ensuremath{\mathrm{TP}_{\mathrm{#1}}}}

\newcommand{\solphaseone}{{\mathcal{S}_1}}
\newcommand{\solphasetwo}{{\mathcal{S}_2}}

\newcommand{\solclasswise}{{\mathcal{S}(\mathcal{A},I_c)}}
\newcommand{\solclassopt}{{\mathcal{S}^\star(I_c)}}

\newtheorem{theorem}{Theorem}[section]
\newtheorem{proposition}[theorem]{Proposition}
\newtheorem{lemma}[theorem]{Lemma}
\newtheorem{corollary}[theorem]{Corollary}
\newtheorem{definition}[theorem]{Definition}
\newtheorem{remark}[theorem]{Remark}
\theoremstyle{definition}

\newcommand{\items}{{\mathcal{I}}}
\newcommand{\classes}{{\mathcal{C}}}
\newcommand{\opt}{OPT}
\newcommand{\UB}{{UB}}
\newcommand{\bincost}{{r}}
\newcommand{\ratio}[1]{R_{#1}}
\newcommand{\ratioinf}[1]{R^{\infty}_{#1}}

\begin{document}
\let\WriteBookmarks\relax
\def\floatpagepagefraction{1}
\def\textpagefraction{.001}

\shorttitle{A Tight 2-Approximation Algorithm for the Bin Packing Problem with Setups}
\shortauthors{Baldacci et al.}

\title[mode=title]{A Tight 2-Approximation Algorithm for the Bin Packing Problem with Setups}

\author[1]{Roberto Baldacci}

\author[2]{Fabio Ciccarelli}
\cormark[1]
\ead{f.ciccarelli@uniroma1.it}

\author[3]{Stefano Coniglio}

\author[2]{Valerio Dose}

\author[2]{Fabio Furini}

\affiliation[1]{organization={College of Science and Engineering, Hamad Bin Khalifa University, Qatar Foundation},
  city={Doha},
  country={Qatar}}

\affiliation[2]{organization={Department of Computer, Control and Management Engineering Antonio Ruberti, Sapienza University of Rome},
  city={Rome},
  country={Italy}}

\affiliation[3]{organization={Department of Economics, University of Bergamo},
  city={Bergamo},
  country={Italy}}

\cortext[1]{Corresponding author}

\begin{abstract}
We study approximation algorithms for the \emph{Bin Packing Problem with Setups} (\BPPS), a generalization of the classical \emph{Bin Packing Problem} (\BPP) in which items are partitioned into classes and activating a class in a bin consumes a setup weight and incurs a setup cost. We show that direct adaptations of Next Fit (NF), First Fit (FF), Best Fit (BF), and Worst Fit (WF), as well as their decreasing-order variants, have unbounded absolute worst-case performance ratios, even with unit-weight items and zero setup costs. 
We then introduce a two-phase algorithm, \TP{\mathcal A}, that packs each class independently with a \BPP{} algorithm \(\mathcal A\) and subsequently merges compatible packing patterns. We prove that the solution returned by \TP{\mathcal A} has cost at most twice the optimum under the assumption that $\mathcal{A}$ produces pairwise merge-maximal solutions, i.e., such that no two packing patterns in the class-wise solution can be feasibly merged. If \(\mathcal A\) also runs in polynomial time, this yields a \(2\)-approximation algorithm for the \BPPS{}.
The factor is tight: the absolute worst-case performance ratio of \TP{\mathcal A} is exactly \(2\), even when \(\mathcal A\) solves every class-wise \BPP{} instance optimally. Since every Any Fit algorithm returns pairwise merge-maximal solutions, it follows that \TP{FF}, \TP{BF}, \TP{WF}, and their decreasing-order variants all have an absolute worst-case performance ratio exactly \(2\).
If, in addition, \(\mathcal A\) is an \(\alpha\)-approximation algorithm with \(\alpha\le2\), we obtain a finer, component-wise guarantee with factor \(2\) for the bin-opening cost and factor \(\alpha\) for the setup-cost component. 
\end{abstract}


\maketitle

\section{Introduction}\label{sec:intro}

Consider an unlimited supply of identical \emph{bins} of capacity \(d \in \mathbb{Z}_{\ge 1}\), each incurring a fixed \emph{bin cost} \(\bincost\in\mathbb{Z}_{\ge 1}\), and a set \(\items = \{1,2,\dots,n\}\) of \(n\) \emph{items}, where each item \(i\in\items\) has weight \(w_i \in \mathbb{Z}_{\ge 1}\).
The items are partitioned into \(m\) (nonempty) \emph{classes}. For each class \(c\in\classes = \{1,2,\dots,m\}\), let \(\items_c\subseteq\items\) denote the set of items of class \(c\). Each class \(c \in \classes\) is characterized by a \emph{setup weight} \(s_c\in\mathbb{Z}_{\ge 0}\) and a \emph{setup cost} \(f_c\in\mathbb{Z}_{\ge 0}\).
If a bin contains at least one item of class \(c\), then class \(c\) is said to be \emph{active} in that bin, the available capacity of the bin is reduced by \(s_c\), and an additional cost \(f_c\) is incurred. A bin is said to be \emph{open} if it contains at least one item; each open bin incurs the fixed bin cost \(\bincost\).

For any subset of items \(S\subseteq\items\), let \(\classes(S):=\{c\in\classes:S\cap\items_c\neq\emptyset\}\) denote the set of classes active in \(S\), and let
\[
\ell(S):=\sum_{i\in S}w_i+\sum_{c\in\classes(S)}s_c
\]
be its \emph{load}. A \textit{packing pattern} is a subset \(S\subseteq\items\) such that \(\ell(S)\le d\). Thus, every packing pattern can be assigned to a single bin. A feasible solution to the \emph{Bin Packing Problem with Setups} (\BPPS{}) is a partition \(\mathcal{S}\subseteq2^{\items}\) of \(\items\) into nonempty packing patterns. Its cost is
\begin{equation}\label{eq:BPPS_solution_cost}
\UB(\mathcal{S})
:=
\sum_{S\in\mathcal{S}}\left(\bincost+\sum_{c\in\classes(S)}f_c\right)
=
\bincost\,|\mathcal{S}|+\sum_{S\in\mathcal{S}}\sum_{c\in\classes(S)}f_c,
\end{equation}
that is, the sum, over all packing patterns, of the cost of the bin to which the pattern is assigned and the setup costs of its active classes. The \BPPS{} consists of finding a feasible solution of minimum cost. We denote by \(\opt(I)\) the optimal value of an instance \(I\) of the \BPPS{}, i.e., the minimum total cost required to pack all items. Hence, for every feasible \BPPS{} solution \(\mathcal{S}\) to \(I\), \(\UB(\mathcal{S})\ge\opt(I)\).
To ensure the feasibility of a \BPPS{} instance, we assume that, for each class \( c \in \classes \), the combined weight of any item \( i \in \items_c \) and its setup weight does not exceed the bin capacity, i.e., \( w_i + s_c \le d \). 
Integrality of the bin cost and of the setup costs is assumed only for ease of exposition: all the results of this paper hold for rational costs as well; only the item weights, the setup weights, and the bin capacity are required to be integer.

The \BPPS{} is strongly \(\textsf{NP}\)-hard, as it admits the classical \textit{Bin Packing Problem} (\BPP{}) as a special case. 
The \BPPS{} has applications in production planning and logistics, and was recently introduced in \cite{BaldacciAndFriends25}---which, to the best of our knowledge, remains the only previous study devoted specifically to this problem. 
The authors of \cite{BaldacciAndFriends25} propose a natural integer linear programming formulation, as well as an arc-flow formulation, for the \BPPS{}, analyze the linear programming relaxations of the two formulations, and derive strengthening valid inequalities.
As far as exact algorithms for the classical \BPP{} are concerned, the current state of the art is represented by branch-price-and-cut methods built on the set-partitioning (or set-covering) reformulation. \citet{WeiLuoBaldacciLim20} solve the pricing problem by an ad hoc label-setting algorithm and strengthen the formulation by means of subset-row inequalities; \citet{BaldacciEtAl24} combine column generation with a pattern-enumeration phase and numerically safe bounding techniques; and \citet{daSilvaSchouery25} propose an extended Ryan--Foster branching scheme together with a fast column generation procedure, improving on the previous best results for the closely related Cutting Stock Problem and for several of its packing relatives.

A problem closely related to the \BPPS{} is the \emph{Knapsack Problem with Setups} (KPS). In the KPS, items are partitioned into classes, and selecting any item of a class incurs both a setup cost and an additional capacity consumption. Thus, although the KPS asks for a maximum-profit subset of items that can be packed together rather than for packing all items into multiple bins, its class-activation structure mirrors that of the \BPPS{}. We refer the reader to \citet{MPV09}, \citet{furini2018exact} and \citet{pferschy2018improved} for formulations, exact algorithms, and approximation results.

Several variants of the \BPP{} are related to the \BPPS{}. In the \emph{Class-Constrained Bin Packing Problem} (CCBPP), items have colors and each bin may contain only a bounded number of colors. \citet{Shachnai2001313} develop approximation schemes for class-constrained packing problems, whereas \citet{ShachnaiTamir04} establish tight online bounds for the unit-size case. \citet{E2006class} study both offline and online algorithms for general item sizes, and \citet{EpsteinImrehLevin10} give an asymptotic fully-polynomial time approximation scheme for a constant number of classes, together with online results and inapproximability bounds. Exact algorithms for the CCBPP are instead developed by \citet{BORGES2020106455}. For the variable-sized version of the CCBPP, \citet{DawandeKalagnanamSethuraman01} proposed an asymptotic polynomial time approximation scheme and two \(3\)-approximation algorithms; \citet{E2006class} provide a second asymptotic polynomial time approximation scheme when the number of classes is bounded by a constant.

Two further \BPP{} variants are close in cost structure to the \BPPS{}. \citet{twigg2015} study a bicriteria problem that asks to minimize both the total number of bins and the number of bins spanned by the items of each color. In the \emph{Bin Packing Problem with Minimum Color Fragmentation}, the number of available bins is instead fixed and the objective is to minimize the total number of times that colors appear in the bins \citep{bergman2019binary,mehrani2022models,BOLOGNA25}. With zero setup weights and unit setup costs, this objective coincides with the class-activation component of the \BPPS{} objective; the \BPPS{} additionally asks to choose how many bins to open and allows class-dependent setup weights and costs.
Lastly, the price-of-clustering literature provides a further connection to the \BPPS{}: it compares an optimal solution obtained by packing clusters of items independently with a globally-optimal packing (see, e.g., \citet{EPSTEIN2021100647,epstein2023}).

Despite the aforementioned streams of work, to the best of our knowledge, no approximation algorithm has been proposed specifically for the \BPPS{}.
Given the rich body of approximation results for the classical \BPP{}---one of the most extensively studied problems in the approximation-algorithms literature, see \cite{coffman2013bin,SpringerRWE_BinPackingOverview}---it is natural to ask whether the \BPPS{} inherits some of its favorable approximability properties. 
In this paper, we take a first exploratory step in this direction and start investigating this question.

The remainder of the paper is organized as follows. Section~\ref{sec:bpp_heuristics} reviews the classical online and offline algorithms for the \BPP{}, together with the performance measures and guarantees used throughout the article; it also reports the immediate consequence that, unless \(\textsf{P}=\textsf{NP}\), no polynomial-time algorithm for the \BPPS{} can have an absolute worst-case performance ratio smaller than \(3/2\). Section~\ref{sec:negative} studies the direct adaptations of these algorithms to the \BPPS{}. We prove that classical \BPP{} algorithms, as well as their decreasing-order variants, have unbounded worst-case performance ratios if straightforwardly applied to the \BPPS{}, even when all items have unit weight and all setup costs are zero. Section~\ref{sec:algorithm} introduces the two-phase algorithm \TP{\mathcal A}, which first packs each class separately by means of a \BPP{} algorithm \(\mathcal A\) and then merges compatible packing patterns. We prove that the solution returned by \TP{\mathcal A} has cost at most twice the optimum whenever \(\mathcal A\) produces pairwise merge-maximal class-wise packings, i.e., solutions in which no two packing patterns can be feasibly merged. We also show that this factor is tight, even if the class-wise subproblems are solved to optimality. It follows that, when \(\mathcal A\) runs in polynomial time, \TP{\mathcal A} is a \(2\)-approximation algorithm for the \BPPS{}.
Finally, we derive a refined guarantee that applies factor \(2\) to the bin-cost component of an optimal solution and the approximation factor of \(\mathcal A\) to its setup-cost component. Section~\ref{sec:conclusions} concludes the paper and discusses directions for further research.

\section{Classical \BPP{} approximation algorithms}\label{sec:bpp_heuristics}

In this section, we review the classical online and offline approximation algorithms for the \BPP{} and recall their performance guarantees. For a detailed review of them and their analyses, we refer the reader to \citet{coffman2013bin,SpringerRWE_BinPackingOverview}. These algorithms provide both the benchmark for the direct \BPPS{} adaptations studied in Section~\ref{sec:negative} and for the algorithms used in Phase~1 of the two-phase algorithm introduced in Section~\ref{sec:algorithm}. Before describing them, we introduce the notation needed to state their approximation properties. In particular, we recall the notions of an \(\alpha\)-approximation algorithm and absolute worst-case performance ratio.

Let \(\mathcal{A}\) be an algorithm for a minimization problem. For any instance \(I\), denote by \(\mathcal{A}(I)\) the value of the solution returned by \(\mathcal{A}\) and by \(\opt(I)\) the optimal value (assumed to be positive). We assume that \(\mathcal{A}\) always returns a feasible solution; hence, for every instance \(I\), it holds that \(\mathcal{A}(I)\ge \opt(I)\) and \(\mathcal{A}(I)/\opt(I)\ge 1\). 
If there exists a finite constant \(\alpha\in[1,+\infty)\) such that
\[
    \frac{\mathcal{A}(I)}{\opt(I)} \le \alpha \quad \text{for all instances } I,
\]
and, in addition, \(\mathcal{A}\) runs in polynomial time, then \(\mathcal{A}\) is called an \textit{\(\alpha\)-approximation algorithm}.

Furthermore, the \textit{absolute worst-case performance ratio} of $\mathcal{A}$ is
\begin{equation} \label{eq:absolute-wcpr}
\ratio{\mathcal A}:=\sup_I \left\{\frac{\mathcal A(I)}{\opt(I)} \right\}\in[1,+\infty].
\end{equation}
Hence, a polynomial time algorithm \(\mathcal A\) is an \(\alpha\)-approximation algorithm if and only if \(\ratio{\mathcal A}\le\alpha\). Nonetheless, note that \(\ratio{\mathcal A}\) is defined for every algorithm \(\mathcal A\), including non-polynomial ones. A factor \(\alpha\) such that \(\mathcal A(I)\le\alpha\,\opt(I)\) for all \(I\) is called
\emph{tight} if it coincides with the ratio in~\eqref{eq:absolute-wcpr}, i.e., if
\(\ratio{\mathcal A}=\alpha\), so that \(\alpha\) cannot be replaced by any smaller value.

\medskip

Before reviewing the classical algorithms, we comment on which performance measure is meaningful for the \BPPS{}. In the \BPP{} literature it is customary to complement the absolute worst-case performance ratio with its asymptotic counterpart, obtained by restricting the supremum in~\eqref{eq:absolute-wcpr} to instances whose optimal value is large. For the \BPPS{} this is vacuous, as the following remark shows.

\begin{remark}\label{rem:cost_scaling}
Call an algorithm for the \BPPS{} \emph{cost-oblivious} if the partition it returns depends on the instance only through \(\items\), \(\classes\), \(d\), \(\boldsymbol{w}\) and \(\boldsymbol{s}\), and not through the bin cost \(\bincost\) or the setup costs \(\boldsymbol{f}\). Every algorithm considered in this paper is cost-oblivious. Given a \BPPS{} instance \(I=(\items,\classes,d,\bincost, \boldsymbol{w},\boldsymbol{s},\boldsymbol{f})\) and \(\lambda\in\mathbb{Z}_{\ge1}\), let \( I(\lambda):=(\items,\classes,d,\lambda\bincost,\boldsymbol{w},\boldsymbol{s},\lambda \boldsymbol{f})\). The two instances have the same feasible solutions, and by~\eqref{eq:BPPS_solution_cost} the cost of each of them is multiplied by \(\lambda\) in \(I(\lambda)\). Hence \(\opt(I(\lambda))=\lambda\,\opt(I)\) and, for a cost-oblivious algorithm \(\mathcal A\), \(\mathcal A( I(\lambda))=\lambda\,\mathcal A(I)\).
Two consequences follow.
\begin{enumerate}
\item[(i)] The ratio \(\mathcal A(I)/\opt(I)\) is invariant under cost scaling, while \(\opt(I(\lambda))\to+\infty\) as \(\lambda\to+\infty\). Restricting the supremum in~\eqref{eq:absolute-wcpr} to instances with large optimal value would therefore return \(\ratio{\mathcal A}\) itself. This is in contrast with the \BPP{}, where the optimal value is the number of open bins and is thus a purely combinatorial quantity, so that letting it grow does constrain the structure of the instance;
\item[(ii)] Additive refinements are impossible: if \(\mathcal A(I)\le\rho\,\opt(I)+ \kappa\) for every instance \(I\), for some constants \(\rho\ge1\) and \(\kappa\ge0\), then applying this inequality to \(I(\lambda)\), dividing by \(\lambda\) and letting \(\lambda\to+\infty\) yields \(\mathcal A(I)\le\rho\,\opt(I)\) for every \(I\). In particular, no \BPPS{} counterpart of bounds such as the one recalled below for FFD, in which an additive term buys a smaller multiplicative factor, can hold for a cost-oblivious algorithm.
\end{enumerate}
\end{remark}

For these reasons, all the results of this paper are stated in terms of the absolute worst-case
performance ratio. A non-degenerate asymptotic notion, obtained by letting the number of open
bins grow, is discussed in Remark~\ref{rem:asymptotic_bins} below.

\medskip

We now briefly review classical \emph{online algorithms} for the \BPP{}, i.e., algorithms that process the items in the given order of the instance.
\begin{itemize}
  \item \textit{Next Fit} (NF) maintains exactly one open bin and places the current item into it if it fits; otherwise, it closes the bin and opens a new one.
  \item \textit{First Fit} (FF) places the current item into the lowest-indexed open bin into which it fits, opening a new bin if necessary.
  \item \textit{Best Fit} (BF) places the current item into an open bin that leaves the minimum residual capacity (equivalently, it chooses a bin with enough residual capacity having maximum load), opening a new bin if necessary.
  \item \textit{Worst Fit} (WF) places the current item into a feasible open bin that leaves the maximum residual capacity (equivalently, it chooses a bin with enough residual capacity having minimum load), opening a new bin if necessary.
\end{itemize}

FF, BF, and WF are \textit{Any Fit} algorithms, i.e., algorithms that open a new bin only when the current item does not fit in any of the previously opened bins; see \citet{Johnson73}. These online algorithms are also commonly used in the \emph{offline} setting by first sorting the items in non-increasing order of weight and then applying the same packing rule.
In this case, they are called \textit{Next Fit Decreasing} (NFD), \textit{First Fit Decreasing} (FFD), \textit{Best Fit Decreasing} (BFD), and \textit{Worst Fit Decreasing} (WFD). Throughout the paper, we assume that items of equal weight are sorted according to their relative order in the input instance. This tie-breaking convention is part of the definition of every decreasing-order algorithm considered here.

All the algorithms mentioned above run in polynomial time and are
\(\alpha\)-approximation algorithms for the \BPP{}; the corresponding values of
\(\alpha\) are recalled below. We refer the reader to the
surveys~\cite{coffman2013bin,SpringerRWE_BinPackingOverview} for further details and
for an overview of other approximation algorithms for the \BPP{}.
In particular, every Any Fit algorithm has an absolute worst-case performance ratio at most \(2\), and both WF and NF have an absolute worst-case performance ratio exactly \(2\); see \citet{BoyarDosaEpstein12} for the exact absolute results and \citet{Johnson1974WorstCasePB} for the classical worst-case analysis. For FF and BF, the tight absolute bounds
\[
   \mathrm{FF}(I)\le \left\lfloor\frac{17\,\opt(I)}{10}\right\rfloor
   \qquad\text{and}\qquad
   \mathrm{BF}(I)\le \left\lfloor\frac{17\,\opt(I)}{10}\right\rfloor
\]
were established by \citet{DosaSgall13} and by \citet{DosaSgall14}, respectively; thus, \(\ratio{\mathrm{FF}}=\ratio{\mathrm{BF}}=17/10\).
For their decreasing variants, we instead have
\[
\ratio{\mathrm{FFD}}=\ratio{\mathrm{BFD}}=\frac32 \, .
\]
These equalities follow from \citet{SimchiLevi1994NewWR}. Moreover, no polynomial-time algorithm for the \BPP{} can achieve an absolute worst-case performance ratio strictly smaller than \(3/2\) unless \(\textsf{P}=\textsf{NP}\); this is a classical consequence of the \(\textsf{NP}\)-completeness of \textsc{Partition}, see \citet{GareyJohnson79}. The tight additive refinement due to \citet{DosaLiHanTuza13} is \(\mathrm{FFD}(I)\le 11/9\,\opt(I)+6/9\).  Since no later result in this paper uses the finer classical bounds for NFD or WFD, we do not report them here.

\medskip

Since the \BPP{} is the special case of the \BPPS{} in which a single class with zero setup weight and zero setup cost is present, the \(3/2\) threshold recalled above transfers verbatim to the \BPPS{}.

\begin{proposition}\label{prop:inapprox}
Unless \(\textsf{P}=\textsf{NP}\), no polynomial-time algorithm for the \BPPS{} has an absolute worst-case performance ratio smaller than \(3/2\).
\end{proposition}

\begin{proof}
Consider any \BPP{} instance with bin capacity \(d\), item set \(\items\) and item weights \(\boldsymbol w\), and let \(I\) be the \BPPS{} instance with the same items and the same bin capacity, a single class \(\classes=\{1\}\) with \(\items_1=\items\), setup weight \(s_1=0\), setup cost \(f_1=0\), and bin cost \(\bincost=1\). Since \(s_1=0\), the standing assumption \(w_i+s_1\le d\) reduces to \(w_i\le d\) for all $i \in \items$, which holds for every feasible \BPP{} instance. The two instances have the same feasible solutions, and \(\UB(\mathcal S)=|\mathcal S|\) for every such solution \(\mathcal S\) by~\eqref{eq:BPPS_solution_cost}. Hence \(\opt(I)\) equals the minimum number of bins required by the original \BPP{} instance. The transformation is computable in linear time, so any polynomial-time algorithm for the \BPPS{} with absolute worst-case performance ratio \(\rho\) yields, by composition, a polynomial-time algorithm for the \BPP{} with the same ratio. The claim thus follows from the \(3/2\) threshold recalled above.
\end{proof}

\section{Poor worst-case performance of classical \BPP{} algorithms on \BPPS{} instances}\label{sec:negative}

In this section, we show that straightforward adaptations of the classical \BPP{} algorithms NF, FF, BF, and WF have no worst-case guarantee on \BPPS{} instances. The same holds for their decreasing-order variants NFD, FFD, BFD, and WFD.
For the \BPPS{}, these algorithms process items in the given order and iteratively apply their packing rule (next, first feasible, best feasible, or worst feasible, respectively), accounting for class setup weights in the residual-capacity calculation: when an item is packed into a bin whose class is not yet active, the corresponding setup cost is incurred, and the setup weight is charged once for that bin, reducing the remaining capacity by the setup weight in addition to the item weight.

To prove that these algorithms have no worst-case guarantee for the \BPPS{}, we propose families of instances that isolate the effect of setup weights by assigning unit weight to every item and setting all setup costs to zero. An optimal solution to these instances consists of one packing pattern for the items of each class; the solutions returned by straightforward extensions of the \BPP{} algorithms to the \BPPS{} instead repeatedly activate the same classes, thereby consuming a noticeable amount of capacity due to the setup weights. Since the number of these class activations grows with the number of items while the optimal packing opens a constant number of bins, the resulting ratio is unbounded. Thanks to employing identical item weights, the same families of instances are directly applicable to the decreasing-order variants.

We start with NF and exhibit a family of \BPPS{} instances on which NF generates one packing pattern per item, whereas an optimal \BPPS{} solution consists of only two packing patterns.

\begin{proposition}\label{prop:NF}
For the \(\mathrm{\BPPS{}}\), the absolute worst-case performance ratio of \(\mathrm{NF}\) is unbounded.
\end{proposition}

\begin{proof}
Fix an even \(n \ge 4\) and consider the \BPPS{} instance \(I_n\) with \(m=2\) classes, \(|\items_1|=|\items_2|=n/2\), and
\[
d=n-1,\qquad w_i=1, \;\; i\in\items,\qquad s_1=s_2=\frac{n}{2}-1,\qquad f_1=f_2=0 .
\]

For each \(c\in\{1,2\}\),
\[
\sum_{i\in\items_c} w_i+s_c=\frac{n}{2}+\Bigl(\frac{n}{2}-1\Bigr)=n-1=d.
\]

Thus, the items of each class form a feasible packing pattern. Moreover, no packing pattern can contain items from both classes, since \(s_1+s_2+2=n>d\). Hence \(\opt(I_n)=2\,\bincost\).

Assume the items are ordered so that their class labels alternate as
$1,2,1,2,\dots$. After packing the first item (class \(1\)), the residual capacity is
\[
d-(1+s_1)=(n-1)-\frac{n}{2}=\frac{n}{2}-1.
\]
The next item (class \(2\)) would require \(1+s_2= n/2>n/2-1\)
units of capacity, so NF opens a new bin. The same argument repeats for every item, and thus \(\mathrm{NF}(I_n)=n\,\bincost\). Therefore,
\[
\frac{\mathrm{NF}(I_n)}{\opt(I_n)}=\frac{n\,\bincost}{2\,\bincost}=\frac{n}{2}\xrightarrow[n\to\infty]{}+\infty.
\]
\end{proof}

An illustration of the difference between an optimal \BPPS{} solution and the solution found by NF for the family of instances in the proof of Proposition~\ref{prop:NF} is given in Figure~\ref{fig:NF}.

\begin{figure}[!htbp]
  \centering
  \begin{center}
\ifx\JPicScale\undefined\def\JPicScale{6}\fi
\unitlength \JPicScale mm
\begin{tikzpicture}[x=21,y=\unitlength,inner sep=0pt]

\def\yMidLow{2.00}   
\def\yMidHigh{3.00}
\def\yMidMid{2.70}

\def\yTopLow{7.40}   
\def\yTopHigh{8.20}
\def\yTopMid{7.80}

\def\ySetLow{2.95}
\def\ySetHigh{3.75}
\def\ySetMid{3.35}

\node[left=0.5cm] at (1,0) {\footnotesize $0$};

\node[left=0.5cm] at (1,1) {\footnotesize $1$};
\node[left=0.5cm] at (1,2) {\footnotesize $2$};
\node[left=0.5cm] at (1,3) {\footnotesize $\frac{n}{2}-2$};
\node[left=0.5cm] at (1,4) {\footnotesize $\frac{n}{2}-1$};
\node[left=0.5cm] at (1,5) {\footnotesize $\frac{n}{2}$};
\node[left=0.5cm] at (1,9) {\footnotesize $n-1$};

\node[below=0.5cm] at (2,0) {bin $1$};
\node[below=0.5cm] at (4,0) {bin $2$};

\node[below=1.2cm] at (3,0) {(a) An optimal \BPPS{} solution};

\draw[fill=white] (1,0) rectangle (3,1);
\node[above right=0.2cm] at (1,0) {\small $w_1$};

\draw[fill=white] (1,1) rectangle (3,2);
\node[above right=0.2cm] at (1,1) {\small $w_3$};

\node at (2,\yMidMid - 0.25) {$\vdots$};

\draw[fill=white] (1,3) rectangle (3,4);
\node[above right=0.2cm] at (1,3) {\small $w_{n-3}$};

\draw[fill=white] (1,4) rectangle (3,5);
\node[above right=0.2cm] at (1,4) {\small $w_{n-1}$};

\draw[fill=white] (1,5) rectangle (3,\yTopLow);
\node[above right=0.2cm] at (1,5) {\small $s_1$};

\draw[fill=white] (1,\yTopHigh) rectangle (3,9);
\draw[fill=white,draw=none] (1,\yTopLow) rectangle (3,\yTopHigh);
\node at (2,\yTopMid) {$///$};

\draw[fill=white] (3,0) rectangle (5,1);
\node[above right=0.2cm] at (3,0) {\small $w_2$};

\draw[fill=white] (3,1) rectangle (5,2);
\node[above right=0.2cm] at (3,1) {\small $w_4$};

\node at (4,\yMidMid - 0.25) {$\vdots$};

\draw[fill=white] (3,3) rectangle (5,4);
\node[above right=0.2cm] at (3,3) {\small $w_{n-2}$};

\draw[fill=white] (3,4) rectangle (5,5);
\node[above right=0.2cm] at (3,4) {\small $w_{n}$};

\draw[fill=white] (3,5) rectangle (5,\yTopLow);
\node[above right=0.2cm] at (3,5) {\small $s_2$};

\draw[fill=white] (3,\yTopHigh) rectangle (5,9);
\draw[fill=white,draw=none] (3,\yTopLow) rectangle (5,\yTopHigh);
\node at (4,\yTopMid) {$///$};

\draw (3,0)--(3,\yMidLow);
\draw (3,\yMidHigh)--(3,\yTopLow);
\draw (3,\yTopHigh)--(3,9);

\def\off{8}

\node[left=0.5cm] at (1+\off,0) {\footnotesize $0$};
\node[left=0.5cm] at (1+\off,1) {\footnotesize $1$};
\node[left=0.5cm] at (1+\off,5) {\footnotesize $\frac{n}{2}$};
\node[left=0.5cm] at (1+\off,9) {\footnotesize $n-1$};

\node[below=1.2cm] at (6+\off,0) {(b)
The NF algorithm solution};

\node[below=0.5cm] at (2+\off,0) {bin $1$};
\node[below=0.5cm] at (4+\off,0) {bin $2$};
\node[below=0.5cm] at (8+\off,0) {bin $n-1$};
\node[below=0.5cm] at (10+\off,0) {bin $n$};

\node at (6+\off,2.3) {$\cdots$};
\node at (6+\off,6.8) {$\cdots$};

\draw[fill=white] (1+\off,0) rectangle (3+\off,1);
\node[above right=0.2cm] at (1+\off,0) {\small $w_1$};

\draw[fill=white] (1+\off,1) rectangle (3+\off,\ySetLow);
\draw[fill=white] (1+\off,\ySetHigh) rectangle (3+\off,5);
\draw[fill=white,draw=none] (1+\off,\ySetLow) rectangle (3+\off,\ySetHigh);
\node at (2+\off,\ySetMid) {$///$};
\node[above right=0.2cm] at (1+\off,1) {\small $s_1$};

\draw[fill=black!20!white] (1+\off,5) rectangle (3+\off,\yTopLow);
\draw[fill=black!20!white] (1+\off,\yTopHigh) rectangle (3+\off,9);
\draw[fill=white,draw=none] (1+\off,\yTopLow) rectangle (3+\off,\yTopHigh);
\node at (2+\off,\yTopMid) {$///$};

\draw[fill=white] (3+\off,0) rectangle (5+\off,1);
\node[above right=0.2cm] at (3+\off,0) {\small $w_2$};

\draw[fill=white] (3+\off,1) rectangle (5+\off,\ySetLow);
\draw[fill=white] (3+\off,\ySetHigh) rectangle (5+\off,5);
\draw[fill=white,draw=none] (3+\off,\ySetLow) rectangle (5+\off,\ySetHigh);
\node at (4+\off,\ySetMid) {$///$};
\node[above right=0.2cm] at (3+\off,1) {\small $s_2$};

\draw[fill=black!20!white] (3+\off,5) rectangle (5+\off,\yTopLow);
\draw[fill=black!20!white] (3+\off,\yTopHigh) rectangle (5+\off,9);
\draw[fill=white,draw=none] (3+\off,\yTopLow) rectangle (5+\off,\yTopHigh);
\node at (4+\off,\yTopMid) {$///$};

\draw[fill=white] (7+\off,0) rectangle (9+\off,1);
\node[above right=0.2cm] at (7+\off,0) {\small $w_{n-1}$};

\draw[fill=white] (7+\off,1) rectangle (9+\off,\ySetLow);
\draw[fill=white] (7+\off,\ySetHigh) rectangle (9+\off,5);
\draw[fill=white,draw=none] (7+\off,\ySetLow) rectangle (9+\off,\ySetHigh);
\node at (8+\off,\ySetMid) {$///$};
\node[above right=0.2cm] at (7+\off,1) {\small $s_1$};

\draw[fill=black!20!white] (7+\off,5) rectangle (9+\off,\yTopLow);
\draw[fill=black!20!white] (7+\off,\yTopHigh) rectangle (9+\off,9);
\draw[fill=white,draw=none] (7+\off,\yTopLow) rectangle (9+\off,\yTopHigh);
\node at (8+\off,\yTopMid) {$///$};

\draw[fill=white] (9+\off,0) rectangle (11+\off,1);
\node[above right=0.2cm] at (9+\off,0) {\small $w_n$};

\draw[fill=white] (9+\off,1) rectangle (11+\off,\ySetLow);
\draw[fill=white] (9+\off,\ySetHigh) rectangle (11+\off,5);
\draw[fill=white,draw=none] (9+\off,\ySetLow) rectangle (11+\off,\ySetHigh);
\node at (10+\off,\ySetMid) {$///$};
\node[above right=0.2cm] at (9+\off,1) {\small $s_2$};

\draw[fill=black!20!white] (9+\off,5) rectangle (11+\off,\yTopLow);
\draw[fill=black!20!white] (9+\off,\yTopHigh) rectangle (11+\off,9);
\draw[fill=white,draw=none] (9+\off,\yTopLow) rectangle (11+\off,\yTopHigh);
\node at (10+\off,\yTopMid) {$///$};

\foreach \x in {3,5,7,9,11} {
  \draw (\x+\off,0)--(\x+\off,\ySetLow);
  \draw (\x+\off,\ySetHigh)--(\x+\off,\yTopLow);
  \draw (\x+\off,\yTopHigh)--(\x+\off,9);
}

\end{tikzpicture}
\end{center}

\caption{Worst-case family of instances for NF used in the proof of Proposition~\ref{prop:NF}. The items are ordered as done in the proof. Blocks represent item weights and the setup weight of their class; the shaded areas indicate unused capacity. Left: an optimal \BPPS{} solution with two packing patterns, one per class, each assigned to a bin.
Right: the NF solution with one packing pattern per item, each assigned to a bin filled up to \(n/2\) by a single item plus its class setup weight.}
  \label{fig:NF}
\end{figure}

We continue with \(\mathrm{FF}\), \(\mathrm{BF}\), and \(\mathrm{WF}\), describing a family of \BPPS{} instances on which these algorithms open a number of bins which grows linearly with the number of items, whereas an optimal \BPPS{} solution consists of only three bins.

\begin{proposition}\label{prop:FF-BF-WF}
For the \(\mathrm{\BPPS{}}\), the absolute worst-case performance ratios of \(\mathrm{FF}\), \(\mathrm{BF}\), and \(\mathrm{WF}\) are unbounded.
\end{proposition}

\begin{proof}
Let \(n\ge 12\) be divisible by \(6\) (the present proof only requires \(n\ge6\); the stronger assumption \(n\ge12\) is used in Remark~\ref{rem:asymptotic_bins}) and consider the \BPPS{} instance \(I_n\) with \(m=3\) classes,
\(|\items_1|=|\items_2|=|\items_3|=n/3\), and
\[
d=\frac{2n}{3}-1,\qquad w_i=1, \;\; i\in\items,\qquad
s_1=s_2=\frac{n}{3}-1,\qquad s_3=\frac{n}{3}-2,\qquad f_c=0, \;\; c\in\classes.
\]

The items of each class form one packing pattern since
\[
\sum_{i\in\items_1} w_i+s_1=\sum_{i\in\items_2} w_i+s_2
=\frac{2n}{3}-1=d,
\qquad
\sum_{i\in\items_3} w_i+s_3=\frac{2n}{3}-2<d.
\]

Hence, \(\opt(I_n)\le 3\,\bincost\). Moreover, no bin can contain items of both classes \(1\) and \(2\), because even packing only one item of each would require
\[
s_1+s_2+2=\Bigl(\frac{n}{3}-1\Bigr)+\Bigl(\frac{n}{3}-1\Bigr)+2=\frac{2n}{3}>d.
\]

Thus, items of classes \(1\) and \(2\) must belong to different
bins; with only two open bins, they would each be full, leaving no capacity for any item of class \(3\). Hence \(\opt(I_n)\ge 3\,\bincost\), and therefore \(\opt(I_n)=3\,\bincost\).

Assume the items are ordered so that their class labels follow the sequence
\((1,2,3,3)\), repeated exactly \(n/6\) times, followed by the remaining \(n/6\) items of class \(1\) and then the remaining \(n/6\) items of class \(2\). When processing each repetition of the sequence, FF, BF and WF open a new bin for the class-\(1\) item and a new bin for the class-\(2\) item since, after packing any of these items, the residual capacity is
\[
d-(1+s_1)=d-(1+s_2)=\Bigl(\frac{2n}{3}-1\Bigr)-\frac{n}{3}=\frac{n}{3}-1,
\]
while an item of the other class would require a capacity of \(1+s_2=1+s_1=\frac{n}{3}>\frac{n}{3}-1\).
Each class-\(3\) item requires \(1+s_3=\frac{n}{3}-1\), so the two class-\(3\) items in the pattern completely fill up the residual space of the two bins just opened; the two bins are indistinguishable for FF, BF and WF, so the argument does not depend on how ties are broken. Hence, when processing each repetition of the sequence, the algorithm opens two new bins. 

After \(n/6\) repetitions of the sequence, \(n/3\) bins have been opened and are full, and all class-\(3\) items are packed. The remaining \(n/6\)
items of class \(1\) fit in one additional bin, and the same holds for those of class \(2\), since
\[
\frac{n}{6}+s_1=
\frac{n}{6}+s_2=
\frac{n}{6}+\Bigl(\frac{n}{3}-1\Bigr)=\frac{n}{2}-1 \le d.
\]

Therefore \(\mathrm{FF}(I_n)=\mathrm{BF}(I_n)=\mathrm{WF}(I_n)=\left(\frac{n}{3}+2\right)\bincost\), and
\[
\frac{\mathrm{FF}(I_n)}{\opt(I_n)}=\frac{\mathrm{BF}(I_n)}{\opt(I_n)}
=\frac{\mathrm{WF}(I_n)}{\opt(I_n)}
=\frac{\left(\frac{n}{3}+2\right)\bincost}{3\bincost}
=\frac{n}{9}+\frac{2}{3}\xrightarrow[n\to\infty]{}+\infty.
\]
\end{proof}

An illustration of the difference between an optimal \BPPS{} solution and the solution found by FF, BF, and WF for the family of instances in the proof of Proposition~\ref{prop:FF-BF-WF} is given in Figure~\ref{fig:FF-BF}.

\begin{figure}[!htbp]
  \centering
  \begin{center}
\ifx\JPicScale\undefined\def\JPicScale{6.5}\fi
\unitlength \JPicScale mm
\begin{tikzpicture}[x=21,y=\unitlength,inner sep=0pt]

\def\xBandEps{0.05}

\def\ySLow{6.20}
\def\ySHigh{7.40}
\def\ySMid{6.80}

\def\ySOneTwoLow{2.20}
\def\ySOneTwoHigh{3.30}
\def\ySOneTwoMid{2.75}

\def\yUpLow{7.00}
\def\yUpHigh{8.00}
\def\yUpMid{7.50}

\def\yHalfItems{2.50}   
\def\yTopTail{6.50}     

\def\yTailLow{3.70}
\def\yTailHigh{4.80}
\def\yTailMid{4.25}

\node[left=0.15cm] at (1,0) {\footnotesize $0$};
\node[left=0.15cm] at (1,1) {\footnotesize $1$};
\node[left=0.15cm] at (1,2) {\footnotesize $2$};
\node[left=0.15cm] at (1,3) {\footnotesize $\frac{n}{3}-2$};
\node[left=0.15cm] at (1,4) {\footnotesize $\frac{n}{3}-1$};
\node[left=0.15cm] at (1,5) {\footnotesize $\frac{n}{3}$};
\node[left=0.15cm] at (1,9) {\footnotesize $\frac{2n}{3}-1$};

\node[below=0.5cm] at (2,0) {bin $1$};
\node[below=0.5cm] at (4,0) {bin $2$};
\node[below=0.5cm] at (6,0) {bin $3$};

\node at (2,2.6) {$\vdots$};
\node at (4,2.6) {$\vdots$};
\node at (6,2.6) {$\vdots$};

\node[below=1.2cm] at (4,0) {(a) An optimal \BPPS{} solution};

\draw[fill=white] (1,0) rectangle (3,1);  \node[above right=0.2cm] at (1,0) {\small $w_{1}$};
\draw[fill=white] (1,1) rectangle (3,2);  \node[above right=0.2cm] at (1,1) {\small $w_{5}$};
\draw[fill=white] (1,3) rectangle (3,4);  \node[above right=0.2cm] at (1,3) {\small $w_{\frac{5n}{6}-1}$};
\draw[fill=white] (1,4) rectangle (3,5);  \node[above right=0.2cm] at (1,4) {\small $w_{\frac{5n}{6}}$};

\draw[fill=white] (1,5) rectangle (3,\ySLow); \node[above right=0.2cm] at (1,5) {\small $s_1$};
\draw[fill=white] (1,\ySHigh) rectangle (3,9);
\draw[fill=white,draw=none] (1-\xBandEps,\ySLow) rectangle (3+\xBandEps,\ySHigh);
\node at (2,\ySMid) {$///$};

\draw[fill=white] (3,0) rectangle (5,1);  \node[above right=0.2cm] at (3,0) {\small $w_{2}$};
\draw[fill=white] (3,1) rectangle (5,2);  \node[above right=0.2cm] at (3,1) {\small $w_{6}$};
\draw[fill=white] (3,3) rectangle (5,4);  \node[above right=0.2cm] at (3,3) {\small $w_{n-1}$};
\draw[fill=white] (3,4) rectangle (5,5);  \node[above right=0.2cm] at (3,4) {\small $w_{n}$};

\draw[fill=white] (3,5) rectangle (5,\ySLow); \node[above right=0.2cm] at (3,5) {\small $s_2$};
\draw[fill=white] (3,\ySHigh) rectangle (5,9);
\draw[fill=white,draw=none] (3-\xBandEps,\ySLow) rectangle (5+\xBandEps,\ySHigh);
\node at (4,\ySMid) {$///$};

\draw[fill=white] (5,0) rectangle (7,1);  \node[above right=0.2cm] at (5,0) {\small $w_{3}$};
\draw[fill=white] (5,1) rectangle (7,2);  \node[above right=0.2cm] at (5,1) {\small $w_{4}$};
\draw[fill=white] (5,3) rectangle (7,4);  \node[above right=0.2cm] at (5,3) {\small $w_{\frac{2n}{3}-1}$};
\draw[fill=white] (5,4) rectangle (7,5);  \node[above right=0.2cm] at (5,4) {\small $w_{\frac{2n}{3}}$};

\draw[fill=white] (5,5) rectangle (7,\ySLow); \node[above right=0.2cm] at (5,5) {\small $s_3$};
\draw[fill=white] (5,\ySHigh) rectangle (7,8);
\draw[fill=white,draw=none] (5-\xBandEps,\ySLow) rectangle (7+\xBandEps,\ySHigh);
\node at (6,\ySMid) {$///$};

\draw[fill=black!20!white] (5,8) rectangle (7,9);

\def\off{8.4}

\def\xA{1}   \def\xB{3}       
\def\xC{3.8} \def\xD{5.8}     

\def\xE{7.2} \def\xF{9.2}     
\def\xG{10.6}\def\xH{12.6}    

\def\xcAB{2}
\def\xcCD{4.8}
\def\xcEF{8.2}
\def\xcGH{11.6}

\def\xDotLeft{3.4}
\def\yDotLow{4.20}
\def\yDotHigh{8.40}

\node[left=0.15cm] at (\xA+\off,0) {\footnotesize $0$};
\node[left=0.15cm] at (\xA+\off,1) {\footnotesize $1$};
\node[left=0.15cm] at (\xA+\off,5) {\footnotesize $\frac{n}{3}$};
\node[left=0.15cm] at (\xA+\off,6) {\footnotesize $\frac{n}{3}+1$};
\node[left=0.15cm] at (\xA+\off,9) {\footnotesize $\frac{2n}{3}-1$};

\node[left=0.15cm] at (\xE+\off,\yHalfItems) {\footnotesize $\frac{n}{6}$};
\node[left=0.15cm] at (\xE+\off,\yTopTail) {\footnotesize $\frac{n}{2}-1$};

\node[below=0.5cm] at (\xcAB+\off,0) {bin $1$};
\node[below=0.5cm] at (\xcCD+\off,0) {bin $\frac{n}{3}$};
\node[below=0.5cm] at (\xcEF+\off,0) {bin $\frac{n}{3}+1$};
\node[below=0.5cm] at (\xcGH+\off,0) {bin $\frac{n}{3}+2$};

\def\xCap{6.8}
\node[below=1.2cm] at (\xCap+\off,0) {(b)
The
FF/BF/WF algorithm solution};

\node at (\xDotLeft+\off,\yDotLow) {$\cdots$};
\node at (\xDotLeft+\off,\yDotHigh) {$\cdots$};
\node at (\xDotLeft+\off,\yDotLow-3.75) {$\cdots$};

\draw[fill=white] (\xA+\off,0) rectangle (\xB+\off,1);
\node[above right=0.2cm] at (\xA+\off,0) {\small $w_{1}$};

\draw[fill=white] (\xA+\off,1) rectangle (\xB+\off,\ySOneTwoLow);
\node[above right=0.2cm] at (\xA+\off,1) {\small $s_{1}$};
\draw[fill=white] (\xA+\off,\ySOneTwoHigh) rectangle (\xB+\off,5);
\draw[fill=white,draw=none] (\xA+\off-\xBandEps,\ySOneTwoLow) rectangle (\xB+\off+\xBandEps,\ySOneTwoHigh);
\node at (\xcAB+\off,\ySOneTwoMid) {$///$};

\draw[fill=white] (\xA+\off,5) rectangle (\xB+\off,6);
\node[above right=0.2cm] at (\xA+\off,5) {\small $w_{3}$};

\draw[fill=white] (\xA+\off,6) rectangle (\xB+\off,\yUpLow);
\node[above right=0.2cm] at (\xA+\off,6) {\small $s_{3}$};
\draw[fill=white] (\xA+\off,\yUpHigh) rectangle (\xB+\off,9);
\draw[fill=white,draw=none] (\xA+\off-\xBandEps,\yUpLow) rectangle (\xB+\off+\xBandEps,\yUpHigh);
\node at (\xcAB+\off,\yUpMid) {$///$};

\draw[fill=white] (\xC+\off,0) rectangle (\xD+\off,1);
\node[above right=0.2cm] at (\xC+\off,0) {\small $w_{\frac{2n}{3}-2}$};

\draw[fill=white] (\xC+\off,1) rectangle (\xD+\off,\ySOneTwoLow);
\node[above right=0.2cm] at (\xC+\off,1) {\small $s_{2}$};
\draw[fill=white] (\xC+\off,\ySOneTwoHigh) rectangle (\xD+\off,5);
\draw[fill=white,draw=none] (\xC+\off-\xBandEps,\ySOneTwoLow) rectangle (\xD+\off+\xBandEps,\ySOneTwoHigh);
\node at (\xcCD+\off,\ySOneTwoMid) {$///$};

\draw[fill=white] (\xC+\off,5) rectangle (\xD+\off,6);
\node[above right=0.2cm] at (\xC+\off,5) {\small $w_{\frac{2n}{3}}$};

\draw[fill=white] (\xC+\off,6) rectangle (\xD+\off,\yUpLow);
\node[above right=0.2cm] at (\xC+\off,6) {\small $s_{3}$};
\draw[fill=white] (\xC+\off,\yUpHigh) rectangle (\xD+\off,9);
\draw[fill=white,draw=none] (\xC+\off-\xBandEps,\yUpLow) rectangle (\xD+\off+\xBandEps,\yUpHigh);
\node at (\xcCD+\off,\yUpMid) {$///$};

\draw[fill=white] (\xE+\off,0) rectangle (\xF+\off,1.00);
\node[above right=0.2cm] at (\xE+\off,0) {\scriptsize $w_{\frac{2n}{3}+1}$};

\draw[fill=white] (\xE+\off,1.65) rectangle (\xF+\off,\yHalfItems);
\node[above right=0.2cm] at (\xE+\off,1.65) {\scriptsize $w_{\frac{5n}{6}}$};

\node at (\xcEF+\off,1.44) {$\vdots$};

\draw[fill=white] (\xE+\off,\yHalfItems) rectangle (\xF+\off,\yTailLow);
\node[above right=0.2cm] at (\xE+\off,\yHalfItems) {\small $s_{1}$};
\draw[fill=white] (\xE+\off,\yTailHigh) rectangle (\xF+\off,\yTopTail);
\draw[fill=white,draw=none] (\xE+\off-\xBandEps,\yTailLow) rectangle (\xF+\off+\xBandEps,\yTailHigh);
\node at (\xcEF+\off,\yTailMid) {$///$};

\draw[fill=black!20!white] (\xE+\off,\yTopTail) rectangle (\xF+\off,\yUpLow);
\draw[fill=black!20!white] (\xE+\off,\yUpHigh) rectangle (\xF+\off,9);
\draw[fill=white,draw=none] (\xE+\off-\xBandEps,\yUpLow) rectangle (\xF+\off+\xBandEps,\yUpHigh);
\node at (\xcEF+\off,\yUpMid) {$///$};

\draw[fill=white] (\xG+\off,0) rectangle (\xH+\off,1.00);
\node[above right=0.2cm] at (\xG+\off,0) {\scriptsize $w_{\frac{5n}{6}+1}$};

\draw[fill=white] (\xG+\off,1.65) rectangle (\xH+\off,\yHalfItems);
\node[above right=0.2cm] at (\xG+\off,1.65) {\scriptsize $w_{n}$};

\node at (\xcGH+\off,1.44) {$\vdots$};

\draw[fill=white] (\xG+\off,\yHalfItems) rectangle (\xH+\off,\yTailLow);
\node[above right=0.2cm] at (\xG+\off,\yHalfItems) {\small $s_{2}$};
\draw[fill=white] (\xG+\off,\yTailHigh) rectangle (\xH+\off,\yTopTail);
\draw[fill=white,draw=none] (\xG+\off-\xBandEps,\yTailLow) rectangle (\xH+\off+\xBandEps,\yTailHigh);
\node at (\xcGH+\off,\yTailMid) {$///$};

\draw[fill=black!20!white] (\xG+\off,\yTopTail) rectangle (\xH+\off,\yUpLow);
\draw[fill=black!20!white] (\xG+\off,\yUpHigh) rectangle (\xH+\off,9);
\draw[fill=white,draw=none] (\xG+\off-\xBandEps,\yUpLow) rectangle (\xH+\off+\xBandEps,\yUpHigh);
\node at (\xcGH+\off,\yUpMid) {$///$};

\end{tikzpicture}
\end{center}

\caption{Worst-case family of instances for FF, BF, and WF used in the proof of Proposition~\ref{prop:FF-BF-WF}. The items are ordered as indicated in the proof.
Left: an optimal \BPPS{} solution has three packing patterns, one per class.
Right: during each of the first $n/6$ repetitions of the $(1, 2, 3, 3)$ sequence, FF, BF, and WF generate one new packing pattern for the class-$1$ item and one for the class-$2$ item. The two class-$3$ items fill up the residual capacities of the associated bins. The remaining $n/6$ items of class $1$ and the remaining $n/6$ items of class $2$ are then packed into one additional bin per class.}
  \label{fig:FF-BF}
\end{figure}

Note that in the instance families considered in the proofs of Propositions~\ref{prop:NF} and~\ref{prop:FF-BF-WF} all items have unit weight. Therefore, sorting them by non-increasing weight does not induce any meaningful order. Under the tie-breaking convention introduced in Section~\ref{sec:bpp_heuristics}, equal-weight items retain the input order used in the two constructions. Hence, the results shown for NF, FF, BF, and WF also hold for the decreasing-order variants NFD, FFD, BFD, and WFD. Thus:

\begin{corollary}\label{cor:decreasing_unbounded}
For the \(\mathrm{\BPPS{}}\), the absolute worst-case performance ratios of \(\mathrm{NFD}\), \(\mathrm{FFD}\), \(\mathrm{BFD}\), and \(\mathrm{WFD}\) are unbounded.
\end{corollary}

\begin{remark}\label{rem:asymptotic_bins}
The families used in Propositions~\ref{prop:NF} and~\ref{prop:FF-BF-WF} have a bounded optimal value, equal to \(2\bincost\) and \(3\bincost\) respectively, so the two statements concern the absolute worst-case performance ratio. By Remark~\ref{rem:cost_scaling}, no stronger conclusion can be drawn from an asymptotic notion based on letting \(\opt(I)\) grow. A non-degenerate
asymptotic measure is instead obtained by letting the \emph{number of open bins} grow, which is the exact counterpart of the classical \BPP{} definition. Let \(\nu(I)\) denote the minimum number of packing patterns over all optimal solutions to a BPPS instance $I$; the minimum is needed because, when setup costs are positive, distinct optimal solutions may use a different number of bins.
Set
\[
\ratioinf{\mathcal A}
:=
\lim_{k\to+\infty}\ \sup_I\left\{\frac{\mathcal A(I)}{\opt(I)}\ :\ \nu(I)\ge k\right\}
\ \le\ \ratio{\mathcal A}.
\]
Propositions~\ref{prop:NF} and~\ref{prop:FF-BF-WF}, being lower bounds on
\(\ratio{\mathcal A}\), do not transfer to \(\ratioinf{\mathcal A}\) by themselves. They do, however, extend by replication. For \(\vartheta\in\mathbb{Z}_{\ge1}\), let \(I_n^{\vartheta}\) be the instance obtained by juxtaposing \(\vartheta\) copies of \(I_n\), each copy having its own classes, and by listing the items copy after copy. When a copy starts, no bin opened while processing the
previous ones has enough residual capacity to activate a class of the new copy: in the family of Proposition~\ref{prop:NF} the only open bin has residual capacity
\(\frac n2-1<\frac n2=1+s_c\), whereas in that of Proposition~\ref{prop:FF-BF-WF} every previously opened bin is either full or has residual capacity \(\frac n6<\frac n3-1\le 1+s_c\), which is where \(n\ge12\) is used. Hence each algorithm processes every copy exactly as it processes \(I_n\), giving \(\mathcal A(I_n^{\vartheta})=\vartheta\,\mathcal A(I_n)\), while \(\opt(I_n^{\vartheta})\le \vartheta\,\opt(I_n)\) because the union of \(\vartheta\) optimal solutions, one per copy, is feasible. Therefore
\[
\frac{\mathcal A(I_n^{\vartheta})}{\opt(I_n^{\vartheta})}\ \ge\ \frac{\mathcal A(I_n)}{\opt(I_n)},
\]
Moreover, \(I_n^{\vartheta}\) contains \(\vartheta n\) unit-weight items and has the same bin capacity \(d\) as \(I_n\), so that every feasible solution to \(I_n^{\vartheta}\) uses at least \(\lceil\vartheta n/d\rceil\) packing patterns; in particular \(\nu(I_n^{\vartheta})\ge\lceil\vartheta n/d\rceil\to+\infty\) as \(\vartheta\to+\infty\). Consequently, the asymptotic worst-case performance ratios of NF, FF, BF, WF and of their decreasing-order variants are unbounded as well.
\end{remark}

\section{Approximation via class-wise optimization and merging}
\label{sec:algorithm}

In this section, we show that a \(2\)-approximation algorithm for the \BPPS{} can be obtained by applying a \BPP{} algorithm \(\mathcal{A}\) \emph{class by class}, followed by a \textit{merging phase} of compatible packing patterns. The analysis relies on the notion of pairwise merge-maximality, which is defined below. We first introduce the associated class-wise \BPP{} instances.

For each class \(c\in\classes\), define the associated \BPP{} instance \(I_c\) with item set \(\items_c\), item weights \(\{w_i\}_{i\in\items_c}\), and bin capacity \(d_c:= d - s_c\). By the standing assumption \(w_i+s_c\le d\), every item of class \(c\) fits into a bin of capacity $d_c$; hence, \(d_c\ge1\) and \(I_c\) is a feasible \BPP{} instance.
A feasible \BPP{} solution for \(I_c\) is a family of nonempty subsets \(\mathcal{S} \subseteq 2^{\items_c}\) forming a partition of \(\items_c\)
(i.e., \(\bigcup_{S\in\mathcal{S}} S = \items_c\),
$S \cap S' = \emptyset$ for all distinct $S,S' \in \mathcal{S}$);
feasibility requires \(\sum_{i\in S} w_i \le d_c\) for all \(S\in\mathcal{S}\).
Given a \BPP{} algorithm \(\mathcal{A}\), we denote by \(\solclasswise\) the solution returned by
\(\mathcal{A}\) on instance \(I_c\), and by \(\overline{\beta}_c:=|\solclasswise|\) its number of open bins. We also denote by \(\solclassopt\) an optimal solution to \(I_c\), and by \(\beta_c:=|\solclassopt|\) the number of bins used by the solution.

\medskip

The following definition formalizes a property, common to many classical \BPP{} approximation algorithms, that is crucial to the theoretical analysis of the proposed algorithm: it characterizes the class-wise packings constructed in Phase~1 of our algorithms and supports all the approximation guarantees established below.

\begin{definition}\label{def:merge_maximal}
Let \(I\) be a \BPP{} instance with bin capacity \(\hat{d}\), and let \(\mathcal S\) be a feasible solution to \(I\). We say that \(\mathcal S\) is \emph{pairwise merge-maximal} if no two of its packing patterns can be feasibly merged, that is, if
\begin{equation}\label{eq:pairwise_nonmergeability}
\sum_{i\in S}w_i+\sum_{i\in S'}w_i>\hat{d}
\qquad
\text{for all distinct }S,S'\in\mathcal S.
\end{equation} 
\end{definition}

Throughout the rest of the paper, we require \(\mathcal A\) to return a pairwise merge-maximal solution to every \BPP{} instance to which it is applied. In particular, when the input is \(I_c\), the packing \(\solclasswise\) returned by \(\mathcal A\) satisfies condition~\eqref{eq:pairwise_nonmergeability} with \(\hat{d} = d - s_c\).

A natural choice for \(\mathcal A\) is an Any Fit algorithm. To see that every packing returned by such an algorithm is pairwise merge-maximal, consider two of its packing patterns \(S_1\) and \(S_2\), with \(S_1\) opened before \(S_2\), and let \(i\) be the first item assigned to \(S_2\). When \(S_2\) is opened, item \(i\) does not fit in \(S_1\); hence, the load of \(S_1\) at that time plus \(w_i\) exceeds the bin capacity. Since the load of \(S_1\) can only increase thereafter and \(i\in S_2\), the final loads of \(S_1\) and \(S_2\) also sum to more than the bin capacity. Thus, \(S_1\) and \(S_2\) cannot be merged. As the choice of the two patterns is arbitrary, every packing returned by an Any Fit algorithm is pairwise merge-maximal.
The same property holds for every optimal \BPP{} packing: if two of its patterns could be feasibly merged, doing so would produce a feasible solution using one fewer bin, contradicting optimality. Finally, pairwise merge-maximality can be enforced on the output of any \BPP{} algorithm by repeatedly merging feasible pairs of patterns until no such pair remains. Since merging never increases the number of open bins, this postprocessing preserves any approximation guarantee enjoyed by the original algorithm, and it runs in polynomial time.

\medskip

We now introduce our two-phase algorithm for the \BPPS{}, denoted by \TP{\mathcal A} and based on a \BPP{} algorithm~\(\mathcal{A}\). As summarized in Algorithm~\ref{alg:TP}, Phase~1 solves the \BPP{} instance associated with each class independently and collects the resulting packing patterns. Phase~2 then repeatedly merges compatible patterns until no further merge is possible.

\medskip

\begin{algorithm}[H]
\caption{The two-phase algorithm \TP{\mathcal A}}\label{alg:TP}
\KwIn{A \BPPS{} instance \(I=(\items,\classes,d,\bincost,\boldsymbol{w},\boldsymbol{s},\boldsymbol{f})\) and a \BPP{} algorithm \(\mathcal A\).}
\KwOut{A feasible \BPPS{} solution \(\mathcal S\).}
\BlankLine
\tcp{Phase 1: class-wise packing}
\(\mathcal S\leftarrow\emptyset\)\;
\ForEach{class \(c\in\classes\)}{
  construct the \BPP{} instance \(I_c\) with items \(\items_c\), weights \(\{w_i\}_{i\in\items_c}\), and capacity \(d-s_c\)\;
  apply \(\mathcal A\) to \(I_c\) and let \(\solclasswise\) be the packing it returns\;
  \(\mathcal S\leftarrow\mathcal S\cup\solclasswise\)\;
}
\(\solphaseone\leftarrow\mathcal S\)\;
\BlankLine
\tcp{Phase 2: merge-maximal merging}
\While{there are distinct \(S,S'\in\mathcal S\) such that \(\ell(S\cup S')\le d\)}{
  choose any such pair \((S,S')\)\;
  \(\mathcal S\leftarrow \left(\mathcal S\setminus\{S,S'\}\right)\cup\{S\cup S'\}\)\;
}
\(\solphasetwo\leftarrow\mathcal S\)\;
\Return \(\solphasetwo\)\;
\end{algorithm}

\medskip

Following Algorithm~\ref{alg:TP}, we denote by \(\solphaseone\) the value of \(\mathcal S\) at the end of Phase~1 and by \(\solphasetwo\) the solution returned at the end of Phase~2. The choice of the pair to merge in Phase~2 is left unspecified: all the theoretical properties proved below hold for any implementation that iterates the merging phase until no pair of packing patterns can be merged. We only require the selection rule not to depend on the bin cost \(\bincost\) nor on the setup costs \(\boldsymbol f\), so that \TP{\mathcal A} is cost-oblivious in the sense of Remark~\ref{rem:cost_scaling}.

\medskip

If \(\mathcal A\) runs in polynomial time, so does \TP{\mathcal A}. Indeed, Phase~1 runs \(\mathcal A\) once per class. Since every Phase-1 pattern is nonempty, \(|\mathcal S_1|\le n\). Phase~2 hence performs at most \(n-1\) merges. At each iteration, all pairs of current patterns can be scanned, and the feasibility of each pair can be tested in polynomial time from their item loads and active-class sets. Thus Phase~2 runs in polynomial time as well.

Both $\solphaseone$ and $\solphasetwo$ are feasible \BPPS{} solutions. In Phase~1, for each class \(c\in\classes\), \TP{\mathcal A} applies \(\mathcal A\) to the \BPP{} instance \(I_c\) with bin capacity \(d - s_c\). Since \(\solclasswise\) is a feasible \BPP{} solution for \(I_c\), its sets form a partition of \(\items_c\).
Because \(\{\items_c\}_{c\in\classes}\) is a partition of \(\items\), we have that \(\solphaseone\) forms a partition of \(\items\).
Moreover, for any \(S\in\solclasswise\), we have \(\sum_{i\in S} w_i \le d-s_c\); hence, \(\ell(S)=\sum_{i\in S} w_i+s_c\le d\).
Therefore \(\solphaseone\) is a feasible \BPPS{} solution. Phase~2 starts from \(\solphaseone\) and repeatedly replaces two packing patterns \(S\) and \(S'\) by \(S\cup S'\) only when \(\ell(S\cup S')\le d\); therefore, by construction, \(\solphasetwo\) is a feasible \BPPS{} solution.
Figure~\ref{fig:TP-FF} illustrates the solutions returned by \TP{\mathcal A} at the end of the two phases on a small \BPPS{} instance when \(\mathcal A=\mathrm{FF}\).

\begin{figure}[!htbp]
  \centering
  \begin{center}
\ifx\JPicScale\undefined\def\JPicScale{4.2}\fi
\unitlength \JPicScale mm
\begin{tikzpicture}[x=21,y=\unitlength,inner sep=0pt,>=stealth]

\def\xA{1.0} \def\xB{3.0}
\def\xC{3.0} \def\xD{5.0}
\def\xE{5.8} \def\xF{7.8}
\def\xG{7.8} \def\xH{9.8}

\node[left=0.5cm] at (\xA,0) {\footnotesize $0$};
\node[left=0.5cm] at (\xA,5) {\footnotesize $5$};
\node[left=0.5cm] at (\xA,10) {\footnotesize $d=10$};
\node at (3.0,10.65) {\footnotesize class $1$};
\node at (7.8,10.65) {\footnotesize class $2$};

\draw[thick,-] (\xA,-0.3)--(\xA,10);
\draw[thick,-] (0.7,0)--(10.3,0);
\foreach \x in {1.0,3.0,5.0,5.8,7.8,9.8} {
  \draw[-] (\x,0)--(\x,-0.3);
}
\foreach \y in {1,...,10} {
  \draw[-] (0.7,\y)--(1.0,\y);
}

\draw[fill=white] (\xA,0) rectangle (\xB,5);
\node[above right=0.07cm] at (\xA,0) {\scriptsize $w_1=5$};
\draw[fill=white] (\xA,5) rectangle (\xB,9);
\node[above right=0.07cm] at (\xA,5) {\scriptsize $w_2=4$};
\draw[fill=white] (\xA,9) rectangle (\xB,10);
\node[above right=0.07cm] at (\xA,9) {\scriptsize $s_1=1$};

\draw[fill=white] (\xC,0) rectangle (\xD,2);
\node[above right=0.07cm] at (\xC,0) {\scriptsize $w_3=2$};
\draw[fill=white] (\xC,2) rectangle (\xD,3);
\node[above right=0.07cm] at (\xC,2) {\scriptsize $s_1=1$};
\draw[fill=black!20!white] (\xC,3) rectangle (\xD,10);

\draw[fill=white] (\xE,0) rectangle (\xF,4);
\node[above right=0.07cm] at (\xE,0) {\scriptsize $w_4=4$};
\draw[fill=white] (\xE,4) rectangle (\xF,7);
\node[above right=0.07cm] at (\xE,4) {\scriptsize $w_5=3$};
\draw[fill=white] (\xE,7) rectangle (\xF,8);
\node[above right=0.07cm] at (\xE,7) {\scriptsize $s_2=1$};
\draw[fill=black!20!white] (\xE,8) rectangle (\xF,10);

\draw[fill=white] (\xG,0) rectangle (\xH,3);
\node[above right=0.07cm] at (\xG,0) {\scriptsize $w_6=3$};
\draw[fill=white] (\xG,3) rectangle (\xH,4);
\node[above right=0.07cm] at (\xG,3) {\scriptsize $s_2=1$};
\draw[fill=black!20!white] (\xG,4) rectangle (\xH,10);

\node[below=0.45cm] at (2.0,0) {bin $1$};
\node[below=0.45cm] at (4.0,0) {bin $2$};
\node[below=0.45cm] at (6.8,0) {bin $3$};
\node[below=0.45cm] at (8.8,0) {bin $4$};
\node[below=1.15cm] at (5.4,0)
  {(a) The Phase~1 solution $\mathcal S_1$};

\def\xI{12.7} \def\xJ{14.7}
\def\xK{14.7} \def\xL{16.7}
\def\xM{16.7} \def\xN{18.7}

\node[left=0.5cm] at (\xI,0) {\footnotesize $0$};
\node[left=0.5cm] at (\xI,5) {\footnotesize $5$};
\node[left=0.5cm] at (\xI,10) {\footnotesize $d=10$};

\draw[thick,-] (\xI,-0.3)--(\xI,10);
\draw[thick,-] (12.4,0)--(19.2,0);
\foreach \x in {12.7,14.7,16.7,18.7} {
  \draw[-] (\x,0)--(\x,-0.3);
}
\foreach \y in {1,...,10} {
  \draw[-] (12.4,\y)--(12.7,\y);
}

\draw[fill=white] (\xI,0) rectangle (\xJ,5);
\node[above right=0.07cm] at (\xI,0) {\scriptsize $w_1=5$};
\draw[fill=white] (\xI,5) rectangle (\xJ,9);
\node[above right=0.07cm] at (\xI,5) {\scriptsize $w_2=4$};
\draw[fill=white] (\xI,9) rectangle (\xJ,10);
\node[above right=0.07cm] at (\xI,9) {\scriptsize $s_1=1$};

\draw[fill=white] (\xK,0) rectangle (\xL,4);
\node[above right=0.07cm] at (\xK,0) {\scriptsize $w_4=4$};
\draw[fill=white] (\xK,4) rectangle (\xL,7);
\node[above right=0.07cm] at (\xK,4) {\scriptsize $w_5=3$};
\draw[fill=white] (\xK,7) rectangle (\xL,8);
\node[above right=0.07cm] at (\xK,7) {\scriptsize $s_2=1$};
\draw[fill=black!20!white] (\xK,8) rectangle (\xL,10);

\draw[fill=white] (\xM,0) rectangle (\xN,2);
\node[above right=0.07cm] at (\xM,0) {\scriptsize $w_3=2$};
\draw[fill=white] (\xM,2) rectangle (\xN,3);
\node[above right=0.07cm] at (\xM,2) {\scriptsize $s_1=1$};
\draw[fill=white] (\xM,3) rectangle (\xN,6);
\node[above right=0.07cm] at (\xM,3) {\scriptsize $w_6=3$};
\draw[fill=white] (\xM,6) rectangle (\xN,7);
\node[above right=0.07cm] at (\xM,6) {\scriptsize $s_2=1$};
\draw[fill=black!20!white] (\xM,7) rectangle (\xN,10);

\node[below=0.45cm] at (13.7,0) {bin $1$};
\node[below=0.45cm] at (15.7,0) {bin $2$};
\node[below=0.45cm] at (17.7,0) {bin $3$};
\node[below=1.15cm] at (15.7,0)
  {(b) The Phase~2 solution $\mathcal S_2$};

\end{tikzpicture}
\end{center}

  \caption{Solution returned by \TP{FF} at the end of the two phases on a \BPPS{} instance with bin capacity \( d = 10 \), number of items \( n = 6 \), number of classes \( m = 2 \) and bin cost $\bincost=10$.  
  The item weights are  $w_1 = 5$, $w_2 = 4$, $w_3 = 2$, $w_4 = 4$, $w_5 = 3$ and $w_6 = 3$, while the setup weights and costs are $s_1 = s_2 = 1$, $f_1 = 2$ and $f_2 = 3$.  
  The items are partitioned into the classes \( \items_1 = \{1,2,3\} \) and \( \items_2 = \{4,5,6\} \). The bins are shown along the horizontal axis, while the vertical axis represents the bin capacity \( d \). For each bin, the figure displays the total weight of the items packed into it and the setup weight of its active classes. Shaded areas indicate unused capacity. In Phase~1 (left), FF is run separately on the two classes and opens four bins, yielding a total cost of $\UB(\solphaseone) = 50$. In Phase~2 (right), the packing patterns containing items $3$ and $6$ are merged. No two of the resulting packing patterns can be further merged, so this is the final solution returned by \TP{FF}, with a total cost of $\UB(\solphasetwo) = 40$.}
  \label{fig:TP-FF}
\end{figure}

We now introduce the notation needed to establish the approximation ratio of \TP{\mathcal A}. Consider an arbitrary feasible \BPPS{} instance $I$, and let \(\mathcal S^\star\) be an optimal solution to it, so that
\(\UB(\mathcal S^\star)=\opt(I)\). For each class
\(c\in\classes\), let
\[
q_c^\star
:=
\bigl| \left\{S\in\mathcal S^\star:c\in\classes(S) \right\}\bigr|
\]
be the number of packing patterns of \(\mathcal S^\star\) in which class \(c\) is
active. Moreover, for any subset of classes \(\tilde{\classes}\subseteq\classes\), define
the load contributed by those classes in \(\mathcal S^\star\) as
\begin{equation}\label{eq:restricted_optimal_load}
\ell^\star(\tilde{\classes})
:=
\sum_{c\in\tilde{\classes}}
\left(\sum_{i\in\items_c}w_i+q_c^\star s_c\right).
\end{equation}
Thus, \(\ell^\star(\tilde{\classes})\) accounts for the weight of all
items belonging to classes in \(\tilde{\classes}\), together with the
setup weight of each such class for every packing pattern in which it is active in
\(\mathcal S^\star\). In particular,
\begin{equation}\label{eq:optimal_total_load}
\ell^\star(\classes)
=
\sum_{S\in\mathcal S^\star}\ell(S)
\le
d\,|\mathcal S^\star|.
\end{equation}

Finally, we call a class \emph{small} if all its items fit into a single bin and \emph{large} otherwise. Accordingly, we partition the set of classes
into the set \(\classes_{\mathrm S}\) of small classes and the set \(\classes_{\mathrm L}\) of large classes, where
\[
\classes_{\mathrm S}
:=
\Bigl\{c\in\classes:\sum_{i\in\items_c}w_i+s_c\le d\Bigr\}
\qquad\text{and}\qquad
\classes_{\mathrm L}
:=
\Bigl\{c\in\classes:\sum_{i\in\items_c}w_i+s_c>d\Bigr\}.
\]

The following lemma, which will be used to prove the approximation ratio of \TP{\mathcal A}, shows that the number of Phase~1 bins used for each large class is less than twice the load contributed by that class in an optimal solution, divided by the bin capacity. This bound follows directly from pairwise merge-maximality of the packing returned by \(\mathcal A\).

\begin{lemma}\label{lem:large_class_bound}
Let \(I\) be a feasible \BPPS{} instance, let \(\mathcal S^\star\) be an optimal solution to \(I\) and, for each class \(c\in\classes\), let \(q_c^\star\) be the number of packing patterns of \(\mathcal S^\star\) in which \(c\) is active. Let \(c\in\classes_{\mathrm L}\) and let \(\solclasswise\) be a pairwise merge-maximal solution to \(I_c\). Then,
\begin{equation}\label{eq:large_class_bound}
\overline{\beta}_c
<
\frac{2}{d}
\left(\sum_{i\in\items_c}w_i+q_c^\star s_c\right).
\end{equation}
\end{lemma}

\begin{proof}

Since \(c\) is a large class, its items do not fit in one bin of capacity \(d-s_c\); hence \(\overline{\beta}_c\ge2\). Let
$
\solclasswise=\{S_1,\ldots,S_{\overline{\beta}_c}\}.
$
Since \(\solclasswise\) is pairwise merge-maximal, by~\eqref{eq:pairwise_nonmergeability} we have
\[
\sum_{i\in S_j}w_i+\sum_{i\in S_h}w_i>d-s_c
\qquad
\text{for all }1\le j<h\le\overline{\beta}_c.
\]
Summing these inequalities over all
\(\binom{\overline{\beta}_c}{2}\) pairs \((S_j,S_h)\), and observing that each
\(S_j\) occurs in exactly \(\overline{\beta}_c-1\) such pairs, yields
\[
(\overline{\beta}_c-1)\sum_{i\in\items_c}w_i
>
\binom{\overline{\beta}_c}{2}(d-s_c),
\]78
and hence \(\sum_{i\in\items_c}w_i>\frac{\overline{\beta}_c}{2}(d-s_c)\). Moreover, the \(q_c^\star\) packing patterns of \(\mathcal S^\star\) in which \(c\) is active induce a feasible packing for \(I_c\), so \(q_c^\star(d-s_c)\ge\sum_{i\in\items_c}w_i\) and therefore \(q_c^\star>\frac{\overline{\beta}_c}{2}\). Consequently,
\begin{align*}
2\left(\sum_{i\in\items_c}w_i+q_c^\star s_c\right)
&=
2\sum_{i\in\items_c}w_i+2q_c^\star s_c\\
&>
\overline{\beta}_c(d-s_c)+\overline{\beta}_c s_c\\
&=
\overline{\beta}_c d,
\end{align*}
which proves~\eqref{eq:large_class_bound}.
\end{proof}

For every small class \(c\in\classes_{\mathrm S}\), pairwise merge-maximality also implies \(\overline{\beta}_c=1\).
If a feasible packing opened at least two bins, the item loads of any two of its bins would sum to at most the total item weight of the class, which is at most \(d-s_c\), contradicting~\eqref{eq:pairwise_nonmergeability}.

Pairwise merge-maximality further implies that every merge performed in Phase~2 involves packing patterns with disjoint sets of active classes. To see this, suppose that two current packing patterns both contain a certain class \(c \in \classes\). Each current packing pattern is the union of some Phase~1 \BPP{} patterns and therefore contains a distinct original pattern of \(\solclasswise\). By~\eqref{eq:pairwise_nonmergeability}, the combined item weight of class \(c\) in these two original patterns exceeds \(d-s_c\). Their union includes those items and incurs setup weight \(s_c\), so its load exceeds \(d\), independently of its other contents. The proposed merge is therefore infeasible. It follows that every feasible Phase~2 merge combines disjoint active-class sets and that the number of activations of each class remains equal to \(\overline{\beta}_c\) after Phase~2.

\medskip

We are now ready to state the main approximation result for \TP{\mathcal A}. The factor \(2\) follows from two related consequences of pairwise merge-maximality. For each large class, Lemma~\ref{lem:large_class_bound} bounds the number of Phase~1 bins in terms of the item and setup load contributed by that class in \(\mathcal S^\star\). For the small classes, after Phase~2 the remaining packing patterns containing only small classes are pairwise non-mergeable; consequently, if at least two such patterns remain, their average load is greater than \(d/2\). These two bounds control the bin-opening cost, while the same class-wise analysis show that no class is activated in the \TP{\mathcal{A}} solution more than twice as often as in \(\mathcal S^\star\), controlling the setup costs.

\begin{theorem}\label{thm:two_approx_anyfit}
Let \(\mathcal A\) be a \BPP{} algorithm that returns a pairwise merge-maximal solution to every input instance. Then, for every feasible \BPPS{} instance \(I\), the solution \(\mathcal S_2\) returned by \TP{\mathcal A} satisfies
\[
   \UB(\mathcal S_2) \le 2\,\opt(I).
\]
\end{theorem}

\begin{proof}
Consider an arbitrary feasible \BPPS{} instance \(I\) and an optimal solution \(\mathcal S^\star\) to it.
We first bound the number of packing patterns in \(\solphasetwo\), by proving that $|\solphasetwo| \le 2 |\mathcal{S}^\star|$.
By Lemma~\ref{lem:large_class_bound},
\begin{equation}\label{eq:large_total_bound}
\sum_{c\in\classes_{\mathrm L}}\overline{\beta}_c
\le
\frac{2}{d}\ell^\star(\classes_{\mathrm L}),
\end{equation}
and the inequality is strict whenever
\(\classes_{\mathrm L}\ne\emptyset\). Moreover,
\begin{equation}\label{eq:lambda_total_load}
\ell^\star(\classes_{\mathrm L})+\ell^\star(\classes_{\mathrm S})
=
\ell^\star(\classes)
\le
d\,|\mathcal S^\star|,
\end{equation}
where the inequality follows from~\eqref{eq:optimal_total_load}.

\smallskip

Let us now consider the partition of \(\solphasetwo\) into
\[
\mathcal F_{\mathrm L}
:=
\{S\in\solphasetwo:\classes(S)\cap\classes_{\mathrm L}\ne\emptyset\}
\]
and
\[
\mathcal F_{\mathrm S}
:=
\solphasetwo\setminus\mathcal F_{\mathrm L}.
\]
Every packing pattern in \(\mathcal F_{\mathrm L}\) contains at least one of the
\(\sum_{c\in\classes_{\mathrm L}}\overline{\beta}_c\) Phase~1 packing patterns dedicated to large classes. Since Phase~2 only merges existing
packing patterns,
\begin{equation}\label{eq:anchored_bins}
|\mathcal F_{\mathrm L}|
\le
\sum_{c\in\classes_{\mathrm L}}\overline{\beta}_c.
\end{equation}
On the other hand, every packing pattern in \(\mathcal F_{\mathrm S}\) contains only small classes.
Moreover, each small class corresponds to exactly one Phase~1 packing pattern and thus occurs in exactly one Phase~2 packing pattern. Hence, the active-class sets of the elements of \(\mathcal F_{\mathrm S}\) are pairwise disjoint subsets of \(\classes_{\mathrm S}\). Consequently,
\begin{align}
\sum_{S\in\mathcal F_{\mathrm S}}\ell(S)
&\le
\sum_{c\in\classes_{\mathrm S}}
\left(\sum_{i\in\items_c}w_i+s_c\right) \notag\\
&\le
\sum_{c\in\classes_{\mathrm S}}
\left(\sum_{i\in\items_c}w_i+q_c^\star s_c\right) \notag\\
&=
\ell^\star(\classes_{\mathrm S}),
\label{eq:small_load_bound}
\end{align}
where the second inequality follows from \(q_c^\star\ge1\) for all $c \in \classes$.

\smallskip

We distinguish two cases. Suppose first that
\(|\mathcal F_{\mathrm S}|\ge2\). At the end of Phase~2, no two distinct packing patterns \(S,S'\in\mathcal F_{\mathrm S}\) can be merged. Moreover, their active class sets are disjoint, because every small class occurs in only one
Phase~2 packing pattern. Thus,
\[
\ell(S)+\ell(S')
=
\ell(S\cup S')
>
d
\qquad
\text{for all distinct }S,S'\in\mathcal F_{\mathrm S}.
\]
Summing these inequalities as in the proof of Lemma~\ref{lem:large_class_bound} over all pairs and observing that each packing pattern appears in exactly \(|\mathcal F_{\mathrm S}|-1\) pairs yields
$
\sum_{S\in\mathcal F_{\mathrm S}}\ell(S)
>
\frac{d}{2} \, |\mathcal F_{\mathrm S}|
$.
Together with~\eqref{eq:small_load_bound}, this gives
\begin{equation}\label{eq:small_number_bound}
|\mathcal F_{\mathrm S}|
<
\frac{2}{d}\ell^\star(\classes_{\mathrm S}).
\end{equation}
Combining~\eqref{eq:large_total_bound},
\eqref{eq:anchored_bins}, and~\eqref{eq:small_number_bound}, we obtain
\begin{align*}
|\solphasetwo|
&=
|\mathcal F_{\mathrm L}|+|\mathcal F_{\mathrm S}|\\
&\le
\sum_{c\in\classes_{\mathrm L}}\overline{\beta}_c
+|\mathcal F_{\mathrm S}|\\
&<
\frac{2}{d}\left(\ell^\star(\classes_{\mathrm L})+\ell^\star(\classes_{\mathrm S})\right)\\
&\le
2|\mathcal S^\star|.
\end{align*}
This proves the desired bound in the first case.

\smallskip

It remains to consider the case \(|\mathcal F_{\mathrm S}|\le1\). If
\(\classes_{\mathrm L}=\emptyset\), then
\(|\mathcal F_{\mathrm L}|=0\) and
\(|\solphasetwo|=|\mathcal F_{\mathrm S}|\le1\), so the desired bound is
immediate. Otherwise, \(\classes_{\mathrm L}\ne\emptyset\) and the
strict version of~\eqref{eq:large_total_bound} gives
\begin{align*}
|\solphasetwo|
&=
|\mathcal F_{\mathrm L}|+|\mathcal F_{\mathrm S}|\\
&\le
\sum_{c\in\classes_{\mathrm L}}\overline{\beta}_c+1\\
&<
\frac{2}{d}\ell^\star(\classes_{\mathrm L})+1\\
&\le
\frac{2}{d}\left(\ell^\star(\classes_{\mathrm L})+\ell^\star(\classes_{\mathrm S})\right)+1\\
&\le
2|\mathcal S^\star|+1.
\end{align*}
Since both \(|\solphasetwo|\) and \(2|\mathcal S^\star|\) are integers,
the strict inequality implies the desired bound in this case as well.
Thus, in all cases,
\begin{equation}\label{eq:number_two_bound}
|\solphasetwo|\le2|\mathcal S^\star|.
\end{equation}

\smallskip
We finally bound the setup-cost contribution. 
Phase~2 only merges packing patterns with disjoint active class sets, and therefore the number of activations of every class in $\solphasetwo$ does not change with respect to the solution obtained after Phase~1. In other words, each class $c \in \classes$ is active exactly in $\overline{\beta}_c$ packing patterns in $\solphasetwo$.
For \(c\in\classes_{\mathrm S}\), \(\overline{\beta}_c=1\le q_c^\star\). For \(c\in\classes_{\mathrm L}\), by
Lemma~\ref{lem:large_class_bound} we have that
\(\overline{\beta}_c<2q_c^\star\). Consequently,
since all setup costs are nonnegative,
\begin{equation}\label{eq:setup_two_bound}
\sum_{S\in\solphasetwo}\sum_{c\in\classes(S)}f_c
=
\sum_{c\in\classes} \overline{\beta}_c \, f_c
\le
2\sum_{c\in\classes}q_c^\star f_c
=
2\sum_{S\in\mathcal S^\star}\sum_{c\in\classes(S)}f_c.
\end{equation}
Combining~\eqref{eq:number_two_bound}
and~\eqref{eq:setup_two_bound}, we conclude that
\begin{align*}
\UB(\solphasetwo)
&=
\bincost\,|\solphasetwo|
+
\sum_{S\in\solphasetwo}\sum_{c\in\classes(S)}f_c\\
&\le
2\bincost\,|\mathcal S^\star|
+
2\sum_{S\in\mathcal S^\star}\sum_{c\in\classes(S)}f_c\\
&=
2\,\UB(\mathcal S^\star) = 2 \, \opt(I) \, .
\end{align*}
\end{proof}

\begin{corollary}\label{cor:two_approx}
    If $\mathcal{A}$ is a polynomial-time algorithm for the \BPP{} that returns a pairwise merge-maximal solution to every input instance, \TP{\mathcal A} is a $2$-approximation algorithm for the \BPPS{}.
\end{corollary}

The factor \(2\) in Theorem~\ref{thm:two_approx_anyfit} is tight. The following construction contains only small classes and does not depend on the Phase~2 merge-selection rule, since the load of every Phase~1 packing pattern is already too large to allow any feasible merge. 
Hence, this guarantee holds for every possible merge-selection rule in Phase~2, provided that merging continues until no pair of packing patterns can be feasibly merged.

\begin{proposition}\label{prop:small_classes_tight}
For every \BPP{} algorithm \(\mathcal A\) that returns pairwise merge-maximal solutions, the absolute worst-case performance ratio of \TP{\mathcal A} is
\[
\ratio{\TP{\mathcal A}}=2.
\]
\end{proposition}

\begin{proof}
For any integer $q \ge 2$, consider a \BPPS{} instance \(I_q\) with bin capacity \(d=2q(q+1)\), bin cost \(\bincost=1\), and \(2q\) classes, each containing $t=q(q+1)$ items. Index the items class by class, so that class \(c \in \classes\) contains the items with indices in \(\{(c-1)t+1,\ldots,ct\}\). Thus, there are \(2qt=2q^2(q+1)\) items in total. Every item $i \in \items$ has weight \(w_i=1\), while every class $c \in \classes$ has setup weight \(s_c=1\) and setup cost \(f_c=0\). Every class is small because
\[
\sum_{i\in\items_c}w_i+s_c=q(q+1)+1\le2q(q+1)=d.
\]
The total item weight of a class is \(q(q+1)=d/2\), so all its items fit in one bin of the associated \BPP{} instance, whose capacity is \(d-s_c=d-1\). Hence, \(\mathcal A\) returns exactly one bin for each class (by pairwise merge-maximality of the class-wise \BPP{} solutions). Phase~1 consequently produces \(2q\) packing patterns, each of load \(\frac{d}{2}+1\). No two of them can be merged in Phase~2, since their union has load
\[
2\left(\frac d2+1\right)=d+2>d.
\]
Thus \(\TP{\mathcal A}(I_q)=2q\), independently of the merge-selection rule.

Notice that the total item weight is
\[
2q\cdot q(q+1)=2q^2(q+1)=qd.
\]
Thus, \(q\) bins would be completely filled by the items alone and could not accommodate any of the setup weights, which are strictly positive. Every feasible solution to $I_q$ therefore contains at least \(q+1\) packing patterns, each of which must be assigned to a bin. This lower bound is attained by opening \(q+1\) bins and, for every \(b\in\{1,2,\ldots,q+1\}\) and \(c\in\{1,2,\ldots,2q\}\), assigning to bin \(b\) the \(q\) class-\(c\) items with indices from \((c-1)t+(b-1)q+1\) through \((c-1)t+bq\). All the \(t=q(q+1)\) items of each class are thereby assigned across the \(q+1\) packing patterns, each having item load \(2q^2\) and setup load \(2q\). Hence every corresponding bin is full, because
\[
2q^2+2q=2q(q+1)=d.
\]
Consequently,
\[
\opt(I_q)=q+1
\qquad\text{and}\qquad
\frac{\TP{\mathcal A}(I_q)}{\opt(I_q)}
=
\frac{2q}{q+1}\xrightarrow[q\to\infty]{}2.
\]
Since the absolute worst-case performance ratio is the supremum over all instances, this family gives \(\ratio{\TP{\mathcal A}}\ge2\). The reverse inequality follows from Theorem~\ref{thm:two_approx_anyfit} and, therefore, \(\ratio{\TP{\mathcal A}}=2\).
\end{proof}

\begin{corollary}\label{corollary:abs_ratio_exact_A}
For every exact \BPP{} algorithm \(\mathcal A\), $\ratio{\TP{\mathcal A}}=2$.   
\end{corollary}

Figure~\ref{fig:TP-small-classes-tight} compares an optimal solution and the solution returned by \TP{\mathcal A} for the family of instances described in the proof of Proposition~\ref{prop:small_classes_tight}.

\begin{figure}[!htbp]
  \centering
  \begin{center}
\ifx\JPicScale\undefined\def\JPicScale{5.0}\fi
\unitlength \JPicScale mm
\begin{tikzpicture}[x=18,y=\unitlength,inner sep=0pt]

\def\xBandEps{0.05}
\def\yWasteLow{7.10}
\def\yWasteHigh{8.20}
\def\yWasteMid{7.65}

\newcommand{\DrawOptimalBin}[6]{%
  \draw[fill=white] (#1,0) rectangle (#2,1);
  \node[above right=0.08cm] at (#1,0) {\tiny $#3$};
  \node at ({(#1+#2)/2},1.50) {$\vdots$};
  \draw[fill=white] (#1,2) rectangle (#2,3);
  \node[above right=0.08cm] at (#1,2) {\tiny $#4$};
  \draw[fill=white] (#1,3) rectangle (#2,4);
  \node[above right=0.08cm] at (#1,3) {\tiny $s_1$};

  \node at ({(#1+#2)/2},5) {$\vdots$};

  \draw[fill=white] (#1,6) rectangle (#2,7);
  \node[above right=0.08cm] at (#1,6) {\tiny $#5$};
  \node at ({(#1+#2)/2},7.50) {$\vdots$};
  \draw[fill=white] (#1,8) rectangle (#2,9);
  \node[above right=0.08cm] at (#1,8) {\tiny $#6$};
  \draw[fill=white] (#1,9) rectangle (#2,10);
  \node[above right=0.08cm] at (#1,9) {\tiny $s_{2q}$};
}

\newcommand{\DrawTPBin}[5]{%
  \draw[fill=white] (#1,0) rectangle (#2,1);
  \node[above right=0.08cm] at (#1,0) {\tiny $#4$};
  \node at ({(#1+#2)/2},2.35) {$\vdots$};
  \draw[fill=white] (#1,3.70) rectangle (#2,4.70);
  \node[above right=0.08cm] at (#1,3.70) {\tiny $#5$};

  \draw[fill=white] (#1,4.70) rectangle (#2,5.40);
  \node[above right=0.08cm] at (#1,4.70) {\tiny $s_{#3}$};

  \draw[fill=black!20!white] (#1,5.40) rectangle (#2,\yWasteLow);
  \draw[fill=black!20!white] (#1,\yWasteHigh) rectangle (#2,10);
  \draw[fill=white,draw=none]
    (#1-\xBandEps,\yWasteLow) rectangle (#2+\xBandEps,\yWasteHigh);
  \node at ({(#1+#2)/2},\yWasteMid) {$///$};
}

\def\xA{0.8} \def\xB{3.3}
\def\xC{4.8} \def\xD{7.3}
\def\xOptDots{4.05}

\node[left=0.15cm] at (\xA,0)  {\footnotesize $0$};
\node[left=0.15cm] at (\xA,10) {\footnotesize $d$};

\DrawOptimalBin{\xA}{\xB}
  {w_1}{w_q}{w_{(2q-1)t+1}}{w_{(2q-1)t+q}}
\DrawOptimalBin{\xC}{\xD}
  {w_{t-q+1}}{w_t}{w_{2qt-q+1}}{w_{2qt}}

\node at (\xOptDots,0.50) {$\cdots$};
\node at (\xOptDots,5.00) {$\cdots$};
\node at (\xOptDots,9.50) {$\cdots$};

\node[below=0.50cm] at (2.05,0) {bin $1$};
\node[below=0.50cm] at (6.05,0) {bin $q+1$};
\node[below=1.20cm] at (4.05,0)
  {(a) An optimal \BPPS{} solution};

\def\xE{11.4} \def\xF{13.9}
\def\xG{13.9} \def\xH{16.4}
\def\xI{18.2} \def\xJ{20.7}
\def\xTPDots{17.3}

\node[left=0.15cm] at (\xE,0)    {\footnotesize $0$};
\node[left=0.15cm] at (\xE,5.40) {\scriptsize $\tfrac{d}{2}+1$};
\node[left=0.15cm] at (\xE,10)   {\footnotesize $d$};

\DrawTPBin{\xE}{\xF}{1}{w_1}{w_t}
\DrawTPBin{\xG}{\xH}{2}{w_{t+1}}{w_{2t}}
\DrawTPBin{\xI}{\xJ}{2q}{w_{(2q-1)t+1}}{w_{2qt}}

\node at (\xTPDots,0.50) {$\cdots$};
\node at (\xTPDots,5.00) {$\cdots$};
\node at (\xTPDots,9.50) {$\cdots$};

\node[below=0.50cm] at (12.65,0) {bin $1$};
\node[below=0.50cm] at (15.15,0) {bin $2$};
\node[below=0.50cm] at (19.45,0) {bin $2q$};
\node[below=1.20cm] at (15.75,0)
  {(b) The \TP{\mathcal A} solution};

\end{tikzpicture}
\end{center}

  \caption{Worst-case family of instances for \TP{\mathcal A} used in the proof of Proposition~\ref{prop:small_classes_tight}. The items are ordered as specified in the proof. Left: an optimal \BPPS{} solution consists of \(q+1\) packing patterns, each containing \(q\) items from every class. Right: \TP{\mathcal A} produces one packing pattern for each of the \(2q\) classes; each pattern contains all \(q(q+1)\) items of one class and has \(d/2-1\) units of unused capacity.}
  \label{fig:TP-small-classes-tight}
\end{figure}
This construction demonstrates that the class-wise decomposition is the main source of loss in \TP{\mathcal A}. Every associated \BPP{} instance is solved to optimality---indeed, all items of a class fit in a single bin---but the Phase~1 packing patterns have load greater than \(d/2\) and Phase~2 cannot merge any pair of them. The global optimum instead splits every class across several packing patterns and exploits the resulting cross-class combinations. Thus, as stated in Corollary~\ref{corollary:abs_ratio_exact_A}, even using an exact Phase~1 algorithm can produce solutions whose value is arbitrarily close to twice the optimum.
Improving the general factor \(2\) therefore requires reconsidering assignments made in Phase~1, rather than just solving the class-wise subproblems more accurately or changing the merge order of Phase~2.

\medskip

Theorem~\ref{thm:two_approx_anyfit} derives the same factor for the bin-opening cost and the setup costs. The factor of the latter, however, can be refined whenever the \BPP{} algorithm used in Phase~1, besides returning pairwise merge-maximal solutions, is an \(\alpha\)-approximation algorithm with \(\alpha\le2\); the refinement is potentially strict when \(\alpha<2\)---as is the case for FF, BF, FFD and BFD---and it is strict whenever, in addition, the setup-cost component of an optimal solution is positive. Indeed, the packing patterns of an optimal \BPPS{} solution in which class \(c \in \classes\) is active induce a feasible solution for \(I_c\), featuring \(q_c^\star\) open bins. Hence \(\beta_c\le q_c^\star\), and an \(\alpha\)-approximation algorithm opens at most \(\alpha q_c^\star\) bins for that class. The factor \(2\) remains necessary for the total number of open bins, but the setup costs can be charged with factor \(\alpha\). This is formalized in the following corollary.

\begin{corollary}\label{cor:refined_cost}
Let \(\mathcal A\) be an \(\alpha\)-approximation algorithm for the \BPP{} with \(\alpha \le 2\) and which returns a pairwise merge-maximal solution to every input \BPP{} instance. Let \(I\) be a feasible \BPPS{} instance, let \(\mathcal S^\star\) be an optimal solution to \(I\) and, for each class \(c\in\classes\), let \(q_c^\star\) be the number of packing patterns of \(\mathcal S^\star\) in which \(c\) is active. Then, the objective value of the solution returned by \TP{\mathcal A} on \(I\) satisfies
\[
\UB(\solphasetwo)
\le
2 \,\bincost \, |\mathcal S^\star| + \alpha \sum_{c\in\classes}q^\star_c \, f_c
=
2 \, \opt(I)-(2-\alpha)\sum_{c\in\classes}q^\star_c \, f_c.
\]
\end{corollary}

\begin{proof}
The bound~\eqref{eq:number_two_bound} gives
\[
\bincost |\solphasetwo|\le 2\bincost |\mathcal S^\star|.
\]
For each class \(c \in \classes\), the \(q_c^\star\) packing patterns of \(\mathcal S^\star\) in which \(c\) is active induce a feasible packing of \(I_c\), and therefore \(\beta_c\le q_c^\star\). Since \(\mathcal A\) is an \(\alpha\)-approximation algorithm for the \BPP{}, it follows that
\[
\overline{\beta}_c\le \alpha\beta_c\le \alpha q_c^\star.
\]
As observed in the proof of Theorem~\ref{thm:two_approx_anyfit}, Phase~2 does not increase the number of activations of any class. It follows that
\[
\sum_{S\in\solphasetwo}\sum_{c\in\classes(S)}f_c
=
\sum_{c\in\classes}\overline{\beta}_c f_c
\le
\alpha\sum_{c\in\classes}q_c^\star f_c.
\]
Adding the two bounds proves the first inequality. The equality follows from
\(
\opt(I)=\bincost |\mathcal S^\star|+\sum_{c\in\classes}q_c^\star f_c
\).
\end{proof}

Note that the refinement vanishes on instances with zero setup costs, which is precisely the case of the tight family of Proposition~\ref{prop:small_classes_tight}. Equivalently, letting \(\sigma(I):=\sum_{c\in\classes}q_c^\star f_c/\opt(I)\in[0,1)\) denote the share of the cost of \(\mathcal S^\star\) due to class activations, Corollary~\ref{cor:refined_cost} reads
\[
\frac{\UB(\solphasetwo)}{\opt(I)}\le 2-(2-\alpha)\,\sigma(I),
\]
indicating that the loss with respect to the optimum decreases as the setup costs become more relevant.

Theorem~\ref{thm:two_approx_anyfit} and Proposition~\ref{prop:small_classes_tight} apply, in particular, to the standard Any Fit algorithms discussed in Section~\ref{sec:bpp_heuristics}. Sorting the items within each class does not affect pairwise merge-maximality: FFD, BFD, and WFD still open a new bin only when the current item fits in none of the previously opened bins. Thus, FF, BF, and WF, as well as their decreasing-order variants, give the following immediate consequence.

\begin{corollary}\label{cor:anyfit_ratios}
The absolute worst-case performance ratios of \TP{FF}, \TP{BF}, \TP{WF}, \TP{FFD}, \TP{BFD}, and \TP{WFD} are equal to \(2\).
\end{corollary}

Finally, the guarantees established in this section hold for the asymptotic measure
of Remark~\ref{rem:asymptotic_bins} as well. The bound of Theorem~\ref{thm:two_approx_anyfit} holds for every instance and, in particular, for instances with arbitrarily many open bins, so that \(\ratioinf{\TP{\mathcal A}}\le2\). Conversely, the family of Proposition~\ref{prop:small_classes_tight} satisfies \(\nu(I_q)=q+1\to+\infty\) as $q$ goes to $+\infty$, so the matching lower bound is asymptotic in nature. Consequently \(\ratioinf{\TP{\mathcal A}}=\ratio{\TP{\mathcal A}}=2\), and the ratios of Corollary~\ref{cor:anyfit_ratios} are asymptotic worst-case performance ratios as well.

\section{Conclusions}\label{sec:conclusions}

We studied approximation algorithms for the \BPPS{}. First, we considered direct extensions of classical \BPP{} algorithms to this variant. Although NF, FF, BF, and WF have finite worst-case performance ratios for the \BPP{}, we proved by counterexample that their straightforward \BPPS{} adaptations, as well as those of their decreasing-order variants, have unbounded absolute worst-case performance ratios. The counterexamples have unit-weight items and zero setup costs, showing that setup weights alone suffice to invalidate the classical guarantees.

We then introduced the two-phase algorithm \TP{\mathcal A}, which packs each class separately in Phase~1 and repeatedly merges compatible packing patterns in Phase~2. We prove that the solution returned by \TP{\mathcal A} has cost at most twice the optimum under the assumption that $\mathcal{A}$ returns pairwise merge-maximal solutions---this property holds for both Any Fit and exact \BPP{} algorithms, and can be obtained by postprocessing any class-wise packing. Hence, if \(\mathcal A\) also runs in polynomial time, \TP{\mathcal{A}} is a \(2\)-approximation algorithm for the \BPPS{}. The bound is tight, even when all classes fit in a single bin, setup weights are equal to one, setup costs are zero, and every class-wise \BPP{} instance is solved to optimality. It follows that \TP{FF}, \TP{BF}, \TP{WF}, and their decreasing-order variants have an absolute worst-case performance ratio exactly \(2\). When \(\mathcal A\) is an \(\alpha\)-approximation algorithm for the \BPP{} with \(\alpha \le 2\), returning pairwise merge-maximal solutions, we further show that the cost of the solution returned by \TP{\mathcal A} is bounded by $2$ times the bin-opening cost of an optimal solution and by $\alpha$ times its setup-cost component. On the negative side, since the \BPP{} is a special case of the \BPPS{}, no polynomial-time algorithm for the latter can have an absolute worst-case performance ratio smaller than \(3/2\) unless \(\textsf{P}=\textsf{NP}\).

Future work includes, first and foremost, understanding whether the factor-\(2\) approximation can be improved. By Proposition~\ref{prop:inapprox}, the absolute approximability threshold of the \BPPS{} lies between \(3/2\) and \(2\), and closing this gap is, in our opinion, the main question left open by this work.
The family of instances we used to prove the tightness of the absolute worst-case performance ratio of \TP{\mathcal A} shows that neither solving the class-wise \BPP{} instances more accurately, even exactly, nor changing the merging rule of Phase~2 would improve the factor 2.
Thus, achieving a better factor likely requires coordinating class activations and item assignments across classes, either within or after Phase~1. 
It would also be natural to investigate stronger guarantees under structural restrictions, such as a bounded number of classes, setup weights or costs being related to bin capacity or bin cost, or additional restrictions on item weights.

\clearpage
\normalsize
\appendix
\section{Notation summary}\label{app:notation}

\noindent Table~\ref{tab:notation} summarizes the notation used throughout the paper.

\footnotesize
\renewcommand{\arraystretch}{1}
\setlength{\tabcolsep}{4pt}
\begin{longtable}[c]{ll}
\caption{\footnotesize{Summary of notation}}\label{tab:notation}\\
\toprule
\textbf{\BPP{} and \BPPS{} input} & \multicolumn{1}{l}{\textbf{Description}} \\
\midrule
\endfirsthead

\caption[]{\footnotesize{Summary of notation (continued)}}\\
\toprule
\textbf{\BPP{} and \BPPS{} input} & \multicolumn{1}{l}{\textbf{Description}} \\
\midrule
\endhead

\bottomrule
\endfoot
\bottomrule
\endlastfoot

$d \in \mathbb{Z}_{\ge 1}$ & Bin capacity \\
$\bincost \in \mathbb{Z}_{\ge 1}$ & Fixed cost of each open bin \\
$n \in \mathbb{Z}_{\ge 1}$ & Number of items \\
$\items=\{1,2,\dots,n\}$ & Set of items \\
$w_i \in \mathbb{Z}_{\ge 1}$ & Weight of item $i\in\items$ \\
$m \in \mathbb{Z}_{\ge 1}$ & Number of classes \\
$\classes=\{1,2,\dots,m\}$ & Set of item classes \\
$\items_c\subseteq\items$ & Set of items belonging to class $c\in\classes$ \\
$s_c \in \mathbb{Z}_{\ge 0}$ & Setup weight of class $c\in\classes$ \\
$f_c \in \mathbb{Z}_{\ge 0}$ & Setup cost of class $c\in\classes$ \\
\midrule
\textbf{\BPPS{} solutions} & \multicolumn{1}{l}{\textbf{Description}} \\
\midrule
$\classes(S)$ & Set of classes active in a subset $S\subseteq\items$ \\
$\ell(S)$ & Load of $S$, including the weights of its items and of its active-class setups \\
$\mathcal{S}\subseteq2^{\items}$ & Feasible \BPPS{} solution, i.e., a partition of $\items$ into nonempty packing patterns \\
$S\in\mathcal{S}$ & Packing pattern assigned to one bin \\
$\UB(\mathcal{S})$ & Cost of a feasible \BPPS{} solution $\mathcal{S}$ \\
$\opt(I)$ & Cost of an optimal solution to a \BPPS{} instance $I$ \\
$\mathcal{S}^\star$ & An optimal solution to the \BPPS{} instance under consideration \\
\midrule
\textbf{Class-wise \BPP{} instances and two-phase heuristic} & \multicolumn{1}{l}{\textbf{Description}} \\
\midrule
$I_c$ & \BPP{} instance induced by class $c$, with item set $\items_c$ and bin capacity $d-s_c$ \\
$d_c:=d-s_c$ & Bin capacity of the class-wise \BPP{} instance $I_c$ \\
$\mathcal{A}$ & \BPP{} algorithm used in Phase~1 of the two-phase heuristic \\
$\solclasswise$ & Solution returned by a \BPP{} algorithm $\mathcal{A}$ to $I_c$ \\
$\overline{\beta}_c:=|\solclasswise|$ & Number of open bins in the class-wise solution $\solclasswise$\\
$\solclassopt$ & Optimal solution to $I_c$ \\
$\beta_c:=|\solclassopt|$ & Optimal number of bins for $I_c$ \\
$\TP{\mathcal{A}}$ & Two-phase \BPPS{} heuristic using algorithm $\mathcal{A}$ in Phase~1\\
$\solphaseone$ & Solution obtained after the class-wise packing phase of $\TP{\mathcal{A}}$ \\
$\solphasetwo$ & Solution returned after the merging phase of $\TP{\mathcal{A}}$ \\
\midrule
\textbf{Approximation analysis} & \multicolumn{1}{l}{\textbf{Description}} \\
\midrule
$I$ & Instance of the \BPPS{} (of the \BPP{}, where explicitly stated) \\
$\mathcal{A}(I)$ & Value of the solution returned by algorithm $\mathcal{A}$ on instance $I$ \\
$\ratio{\mathcal{A}}$ & Absolute worst-case performance ratio of an approximation algorithm $\mathcal{A}$ \\
$\ratioinf{\mathcal{A}}$ & Asymptotic worst-case performance ratio of an approximation algorithm $\mathcal{A}$ \\
\(\nu(I)\) & Minimum number of packing patterns over all optimal solutions to a
\BPPS{} instance \(I\) \\
$\alpha$ & Approximation ratio of an approximation algorithm $\mathcal{A}$ \\
$q_c^\star$ & Number of packing patterns of an optimal \BPPS{} solution $\mathcal{S}^\star$ in which class $c$ is active \\
$\ell^\star(\tilde{\classes})$ & Load contributed in an optimal \BPPS{} solution $\mathcal{S}^\star$ by the classes in $\tilde{\classes}\subseteq\classes$ \\
$\classes_{\mathrm S}$ & Set of small classes, i.e., those whose items fit in a single bin \\
$\classes_{\mathrm L}$ & Set of large classes, i.e., those whose items do not fit in a single bin \\
$\mathcal{F}_{\mathrm S}$ & Packing patterns of $\solphasetwo$ containing items of small classes only \\
$\mathcal{F}_{\mathrm L}$ & Packing patterns of $\solphasetwo$ containing items of at least one large class \\
\end{longtable}
\normalsize

\end{document}